\documentclass[review]{elsarticle}
\usepackage{amsmath,amssymb,amsthm}
\usepackage[mathscr]{eucal}
\usepackage{rotate,graphics,epsfig}
\usepackage{color}
\usepackage{hyperref}

\newtheorem{The}{Theorem}
\newtheorem{Cor}[The]{Corollary}
\newtheorem{Pro}[The]{Proposition}

\theoremstyle{definition}
\newtheorem{Rem}{Remark}
\numberwithin{Rem}{section}
\newtheorem{Exa}{Example}
\numberwithin{Exa}{section}
\newtheorem{Def}{Definition}
\numberwithin{Def}{section}
\numberwithin{equation}{section}
\numberwithin{The}{section}

\newcommand{\be}{\begin{eqnarray}}
\newcommand{\ee}{\end{eqnarray}}

\newcommand{\by}{\begin{eqnarray*}}
\newcommand{\ey}{\end{eqnarray*}}

\newcommand{\bn}{\begin{enumerate}}
\newcommand{\en}{\end{enumerate}}

\newcommand{\bi}{\begin{itemize}}
\newcommand{\ei}{\end{itemize}}

\def\oo{\infty}

\def\Fbar{{\overline F}}

\def\frac#1#2{{#1 \over #2}}

\def \P {\mathbb{P}}

\def\xx{\boldsymbol{x}}
\def\XX{\boldsymbol{X}}
\def\mmu{\boldsymbol{\mu}}
\def\ddelta{\boldsymbol{\delta}}

\def\ttheta{\boldsymbol{\theta}}
\def\yy{\boldsymbol{y}}
\def\zz{\boldsymbol{z}}
\def\YY{\boldsymbol{Y}}

\def\ww{\boldsymbol{w}}

\def\Zero{\boldsymbol{0}}

\begin{document}

\begin{frontmatter}
\title{Tail Densities of Skew-Elliptical Distributions}
\author[label1]{Harry Joe}
\ead{Harry.Joe@ubc.ca}
\author[label2]{Haijun Li\corref{cor2}}
\ead{lih@math.wsu.edu}
\cortext[cor2]{Corresponding author}

\address[label1]{Department of Statistics,
University of British Columbia, Vancouver, BC, V6T1Z4, Canada}
\address[label2]{Department of Mathematics and Statistics, Washington State
University, Everett, WA 98201, U.S.A.}


\begin{abstract}
Skew-elliptical distributions constitute a large class of multivariate distributions that account for both skewness and a variety of tail properties. 
{This class has simpler representations in terms of densities rather
than cumulative distribution functions, and the tail density approach has 
previously been developed to study tail properties when multivariate
densities have more tractable forms. The special skew-elliptical structure 
allows for derivations of specific forms for}
the tail densities for those skew-elliptical copulas
that admit probability density functions, under heavy and light tail conditions on density generators. The tail densities of skew-elliptical copulas are explicit and depend only on tail properties of the underlying density generator and  conditions on the skewness parameters. 
In the heavy tail case 
skewness parameters affect tail densities of the skew-elliptical copulas more profoundly than that in the light tail case, 
whereas in the latter case the tail densities of skew-elliptical copulas are only proportional to the tail densities of symmetrical elliptical copulas. Various examples, including tail densities of skew-normal and skew-t distributions, are given. 
\end{abstract}


\end{frontmatter}


\section{Introduction}
\label{S1}

Several classes of multivariate skew-elliptical distributions have been
introduced and studied in the literature, see, e.g.,  \cite{AD96, BD2001,
AC03, Genton05, Genton06}. These multivariate distributions can account for
both skewness and heavy tails, and have found widespread application in
various areas \cite{Azzalini05, Genton04, Azzalini2013, Yoshiba2018}. One general class of skew distributions consists of selection distributions, which are motivated from various applications in modeling multivariate observations under constraints.

It is our aim 
{to use}
the density approach to establish an extreme value theory for the skew-elliptical distributions introduced in \cite{BD2001} for selection models. We derive the tail density, i.e., the decay rate of some multivariate skew-elliptical  distributions that admit probability density functions, under the conditions that (1) density generators are  regularly varying (the heavy tail case) or (2) density generators are in the quadratic max-domain of attraction for the Gumbel distribution (the light tail case). By using the copula approach, we then drive the scale-invariant upper and lower tail densities, i.e., the decay rates of skew-elliptical copula densities. The scale-invariant tail behaviors can be directly derived 
for copula families that have tractable forms. However for distributions such as skew-elliptical distributions, cumulative distributions and the corresponding copulas may not have closed forms. In contrast, our density-based approach becomes particularly useful in analyzing tail behaviors of skew-elliptical distributions {for which} only joint densities are tractably available.

Our results show that skewness affects tail densities of the skew-elliptical copulas more profoundly in the heavy tail case due to slow tail decays, whereas skew tail densities in the light tail case are proportional to the symmetric tail densities, as are illustrated in skew-normal distributions. 
The logarithms of extreme value distributions and of extreme value copulas can be expressed as integrals of tail densities, and respectively, integrals of copula tail densities, and thus our results provide explicitly distributions for multivariate extremes of observations that are modeled by skew-elliptical densities. Various tail dependence parameters
are obtained for skew-elliptical distributions as by-products.

{The remainder of this} paper is organized as follows. The preliminaries on skew-elliptical distributions and terminology for tail densities of copulas are presented in Section \ref{S15}. The tail densities of the skew-elliptical copulas are derived in Section \ref{sec-theory} under heavy and light tail assumptions on density generators. The results are illustrated for skew-normal and skew-t distributions and also extended to a finite mixture case in Section \ref{sec-examples}. All the proofs are detailed in the Appendix.

\section{Preliminaries: Skew-Elliptical Distributions and Copula Tail Densities}
\label{S15}

In this section, we introduce the notation and {background} for the main results,
and review {some relevant} properties for elliptical
and skew-elliptical distributions. We also list {some relevant} results obtained in \cite{LH14} for the relations of tail densities between multivariate distributions and their copulas. 
{The subsections have topics of skew-elliptical distributions, copulas and tail properties, 
and tail density}.

In what follows, 
vectors are row vectors (typically $1\times d$) shown as boldfaced letters
with elements that are not boldfaced, and so transposes become column vectors. 
Random variables/vectors are upper-case letters and arguments of functions
are lower-case letters. 
Two functions $f,g: {\mathbb{R}}\to \mathbb{R}$ are said to be tail equivalent, denoted by  $f(x)\sim g(x)$ as $x\to a$, $a\in \overline{\mathbb{R}}=\mathbb{R}\cup \{+\infty\}$, if $\lim_{x\to a}[f(x)/g(x)]= 1$. A univariate Borel-measurable function $V: \mathbb{R}_+\to \mathbb{R}_+$ is said to be regularly varying (RV) at $\infty$ with tail index $\rho\in \mathbb{R}$, denoted by $V\in \mbox{RV}_\rho$, if $V(tx)/V(t)\to x^\rho$ as $t\to \infty$ for any $x>0$. See, for example,  \cite{Resnick07} for details on these notions and on theory of regular variation. 

\subsection{Skew-elliptical distributions}

A general class of skew-elliptical distributions has been  constructed in \cite{BD2001} using conditional distributions  from symmetrical elliptical distributions.   A $d$-dimensional random vector ${  \XX}$ is said to be {elliptically distributed}  with mean vector ${  \mmu}\in \mathbb{R}^d$ and dispersion matrix $\Sigma\in \mathbb{R}^{d\times d}$ that is nonnegative-definite, if the characteristic function $\varphi_{{  {\XX-\mmu}}}$ of $  {\XX-\mmu}$ is a function of the quadratic form ${ \zz}\Sigma{\,  \zz}^\top$, where $\zz^\top$ denotes the transpose of a row vector $\zz\in \mathbb{R}^d$.
That is, $\varphi_{  {\XX-\mmu}}({ \zz})=\phi({ \zz}\Sigma{\, \zz}^\top)$, $\zz\in \mathbb{R}^d$,  for a Borel-measurable function $\phi: [0,\infty)\to \mathbb{R}$, which is called a  characteristic generator. We assume that the density of $\XX$ exists and is given by
\begin{equation}
f_{\XX}(\xx; \mmu, \Sigma)= |\Sigma|^{-1/2}g_{d}((\xx-\mmu)\Sigma^{-1}(\xx-\mmu)^\top),~~\xx\in \mathbb{R}^d,
\end{equation}
where $g_{d}(\cdot)$ is known as the density generator that satisfies that
\begin{equation}
\int_0^\infty r^{d/2-1}g_{d}(r){\rm d}r = \frac{\Gamma(d/2)}{\pi^{d/2}}. 
\label{generator integral}
\end{equation}
 The elliptically distributed random vectors with densities  are denoted as ${  \XX}\stackrel{D}{\sim} E_d(\mmu, \Sigma, g_{d})$. A comprehensive review of  characterizations and  properties of elliptical distributions can be found in \cite{FKN90}.  We assume throughout that the density generator $g_{d}(t)$ is continuous over $(0, \infty)$ and eventually non-increasing as $t\to \infty$. Note that this assumption is rather weak and all commonly used density generators are continuous functions that eventually decreasing to zero.

For  skew-elliptical distributions, consider a $(d+1)$-dimensional random vector $(X_0, \XX)\stackrel{D}{\sim} E_{d+1}({  \mmu}^*, \Sigma^*, g_{d+1})$, where $\XX=(X_1, \dots, X_d)$,
 $\mmu^*=(0, \mmu)$, $\mmu=(\mu_1, \dots, \mu_d)$,  $g_{d+1}(\cdot)$ is a density generator, and the dispersion matrix $\Sigma^*$ has the form
 \[\Sigma^*=
 \begin{pmatrix}
  1 & \ddelta \\
  \ddelta^\top & \Sigma
  \end{pmatrix}
 \]
 with $\ddelta = (\delta_1, \dots, \delta_d)\in \mathbb{R}^d$ such that $\Sigma^*$ is positive-definite. The distribution of $\YY:=[\XX|X_0>0]$ is called in \cite{BD2001} a skew-elliptical distribution and denoted as $\YY\stackrel{D}{\sim} SE_{d}({  \mmu}, \Sigma, g_{d+1}, \ddelta)$, where $\ddelta$ is known as the vector of skewness parameters. Note that marginally, $\XX\stackrel{D}{\sim} E_d({  \mmu}, \Sigma, g_{d})$, where 
 \begin{equation}
 \label{g}
 g_{d}(s)= \frac{\pi^{1/2}}{\Gamma(1/2)}\int_s^\infty (r-s)^{-1/2}g_{d+1}(r){\rm d}r=
 2\int_0^\infty g_{d+1}(r^2+s){\rm d}r, ~~s\ge 0,
 \end{equation}
 due to the fact $\Gamma(1/2)=\pi^{1/2}$; 
see (2.23) in \cite{FKN90}.
Conditioning on $\XX=\yy$, the conditional variable $[X_0|\XX=\yy]$ has a univariate elliptical distribution. It follows from Theorem 2.18 of \cite{FKN90} that  $\mathbb{P}(X_0>0|\XX=\yy)=F((\yy-\mmu)\ttheta^\top; g_{q(\yy)})$, where $F(\cdot\ ; g_{q(\yy)})$ is the cumulative distribution function of the univariate elliptical distribution $E_1(0, 1, g_{q(\yy)})$, and $q(\yy)=(\yy-\mmu)\Sigma^{-1}(\yy-\mmu)^\top$, and 
\begin{equation}
\label{eq-ttheta}
g_{q(\yy)}(s)=\frac{g_{d+1}(s+q(\yy))}{g_{d}(q(\yy))}, ~s\ge 0,~~\mbox{and}~\ttheta=\frac{\ddelta\Sigma^{-1}}{(1-\ddelta \Sigma^{-1}\ddelta^\top)^{1/2}}.
\end{equation}
It then follows that the density of $\YY$ can be written as
\begin{equation}
\label{density Y}
f_{\YY}(\yy)=2f_{\XX}(\yy; \mmu, \Sigma)F\big((\yy-\mmu)\ttheta^\top; g_{q(\yy)}\big), ~\yy\in \mathbb{R}^d. 
\end{equation}
As in \cite{BD2001}, we assume that $\ddelta \Sigma^{-1}\ddelta^\top<1$ (equivalently, the positive definiteness of the $(d+1)\times(d+1)$
dispersion matrix $\Sigma^*$), and thus $\ttheta$ in \eqref{eq-ttheta} is well-defined. If $\ddelta = \Zero$, then  $\ttheta=\Zero$ and thus $f_{\YY}(\cdot)=f_{\XX}(\cdot; \mmu, \Sigma)$ in this case. In general, $f_{\YY}(\cdot)$ is skewed. It follows from \eqref{eq-ttheta} and \eqref{density Y} that the density $f_{\YY}(\cdot)$ can be expressed more explicitly as follows:
\begin{equation}
\label{density Y explicit}
f_{\YY}(\yy)=2|\Sigma|^{-1/2}\int_{-\infty}^{(\yy-\mmu)\ttheta^\top}g_{d+1}(r^2+(\yy-\mmu)\Sigma^{-1}(\yy-\mmu)^\top){\rm d}r. 
\end{equation}
The examples of \eqref{density Y} include skew-normal and skew-t distributions introduced and studied in \cite{AD96, AC03}. 
The distributional properties of this class of  skew-elliptical distributions can be found in detail in \cite{BD2001}. 


\subsection{Copulas and tail properties}

The tail properties, such as scale-invariant tail dependence and tail asymmetry,  can be derived from joint tails of the underlying copulas.
Let $F$ denote a $d$-dimensional distribution function with continuous marginal distributions $F_1, \dots, F_d$ and corresponding survival functions $\Fbar_1,\ldots,\Fbar_d$.  Formally, a copula $C$ is a multivariate distribution with standard uniformly distributed margins on $[0,1]$. Sklar's theorem (see, e.g., Section 1.3 in \cite{Joe2014}) states that every {continuous} multivariate distribution $F$ with margins $F_1,\dots, F_d$ can be written as $F(x_1,\dots, x_d)=C(F_1(x_1),\dots,F_d(x_d))$ for some $d$-dimensional copula $C$,
{where}
 $$C(u_1,\dots,u_d)={F}(F_1^{-1}(u_1), \dots, F_d^{-1}(u_d))$$
{and} $F_i^{-1}(u_i)$ is the quantile function of the $i$-th margin
for $1\le i\le d$. 
Let $(U_1,\dots,U_d)$ denote a random vector with distribution $C$ and $U_i, 1\le i\le d$, being uniformly distributed on $[0,1]$. The survival copula $\widehat{C}$ is defined as follows:
\begin{equation}
\label{dual}
\widehat{C}(u_1,\dots,u_n)=\P(1-U_1\le u_1,\dots,1-U_n\le u_n)=\overline{C}(1-u_1,\dots,1-u_n)
\end{equation}
where $\overline{C}$ is the joint survival function of $C$. The survival copula $\widehat{C}$ can be used to transform lower tail properties of $(U_1,\dots,U_d)$ into the corresponding upper tail properties of $(1-U_1,\dots,1-U_d)$. 
{We} assume throughout that the density function $f$ of $F$ exists and 
the density $c(\cdot)$ of copula $C$ exists, and that $c(\cdot)$ is continuous in small open neighborhoods near and not containing
the corner points $\boldsymbol 0$ and ${\boldsymbol 1}=(1, \dots, 1)$.

The upper {\em  tail density of $C$ with tail order $\kappa_U$}, introduced in \cite{LH14},
 is defined as follows:
 \begin{equation}
 \lambda_U({\boldsymbol w}; \kappa_U):= \lim_{u\to 0^+}
\frac{c(1-u{\boldsymbol w})}{u^{\kappa_U-d}\ell_U(u)}>0, ~ {\boldsymbol w}=(w_1,\dots,w_d)\in [0,\infty)^d\backslash \{\boldsymbol 0\},\label{U tail}
\end{equation}
provided that the non-zero limit exists for some $\kappa_U\ge 1$ and some function ${\ell_U}(\cdot)$ that is slowly varying at $0$ (i.e., $\ell_U(us)/\ell_U(u)\to 1$ for any $s>0$ as $u\to 0$). Similarly, the lower {tail density of $C$ with tail order $\kappa_L$} 
 is defined as follows:
 \begin{equation}
 \lambda_L({\boldsymbol w}; \kappa_L):= \lim_{u\to 0^+}
\frac{c(u{\boldsymbol w})}{u^{\kappa_L-d}{\ell_L}(u)}>0, ~ {\boldsymbol w}=(w_1,\dots,w_d)\in [0,\infty)^d\backslash \{\boldsymbol 0\},\label{L tail}
\end{equation}
provided that the limit exists for some $\kappa_L\ge 1$ and some slowing varying function $\ell_L(\cdot)$. Clearly, the tail density functions are homogeneous; that is,
\begin{equation}
\label{scaling}
\lambda_U(t{\boldsymbol w}; \kappa_U)=t^{\kappa_U-d}\lambda_U({\boldsymbol w}; \kappa_U), ~~\lambda_L(t{\boldsymbol w}; \kappa_L)=t^{\kappa_L-d}\lambda_L({\boldsymbol w}; \kappa_U), 
\end{equation}
for any $t>0$ and ${\boldsymbol w}\in [0,\infty)^d\backslash\{\boldsymbol 0\}$. The tail densities \eqref{U tail} and \eqref{L tail} describe, respectively, the scale-invariant decay rates of a copula density near the upper corner $\bf 1$ and the lower corner $\bf 0$. 
The copula tail densities are a key ingredient  in tail risk analysis since the tail risk measures are often expressed in terms of tail densities of the multivariate copulas of underlying loss distributions \cite{BE07, JLN10, JL11, ZL12}.  
The slowly varying function $\ell(\cdot)$ is a non-zero constant in most of the tail dependence and tail risk measures used in practice (e.g., $\ell(u)=1$ in \cite{Yoshiba2018}), and in contrast a general slowly varying function $\ell(\cdot)$ is employed in \eqref{U tail} and \eqref{L tail} for additional flexibility. 
Since the lower tail density of a copula $C$ is the upper tail density of the survival copula $\widehat{C}$ (see \eqref{dual}), we focus only on the upper tail density. The properties of copula tail densities and their relations to tail densities of multivariate extreme value distributions can be found in \cite{LW2013, HJ10, HJL12, LH14}.

\subsection{Tail density results for the Fr\'echet and Gumbel cases}

{The tail densities have special forms when}
the univariate marginal distributions $F_i$'s belong to the max-domain of attraction of one of the following distribution families: Fr\'echet, Gumbel or Weibull distribution.
We only focus on Fr\'echet and Gumbel cases, because the Weibull case is similar to the Fr\'echet case.  We also assume  that for the upper tail case, margins $F_1, \dots, F_d$ are  right-tail equivalent in the sense that
\begin{equation}
	\label{equiv}
	f_i(t)\sim a_if_1(t),~\mbox{as}~t\to \infty,~~1\le i\le d, 
\end{equation}
where $f_i$ is the density of the $i$-th margin $F_i$, and $0<a_i<\infty$ is a constant,  $1\le i\le d$,  with $a_1=1$. Note that \eqref{equiv} implies that the usual tail equivalence $\Fbar_i(t)\sim a_i\Fbar_1(t)$, $1\le i\le d$, as $t\to \infty$, but the reverse may not be true because some densities may have dampened oscillations as $t\to \infty$. The following two results, obtained in \cite{LH14}, specify the relations among the tail densities of multivariate distributions and their copulas. 
{They are useful for establishing joint tail properties if the joint
density has a tractable form but not the joint cumulative distribution
function.}

\begin{Pro}
	\label{F} 
	\rm  (Fr\'echet case) 
	Assume that marginal densities are right-tail equivalent in the sense of \eqref{equiv}. 
	If the marginal density $f_i\in \mbox{RV}_{-\gamma-1}$, $1\le i\le d$, $\gamma>0$, and the limit \eqref{U tail} with tail order $\kappa_U$ holds locally uniformly for all ${\boldsymbol w}\in [0,\infty)^d\backslash \{\boldsymbol 0\}$, then the tail density $\lambda(\cdot)$ of $F$ can be expressed in terms of the upper tail density $\lambda_U(\cdot)$ of $C$ as follows: for any ${\boldsymbol w}=(w_1, \dots, w_d)>{\boldsymbol 0}$, 
	\begin{eqnarray}
		\lambda({\boldsymbol w}) &:=& \lim_{t\to \infty}\frac{f(t{\boldsymbol w})}{t^{-d}V^{\kappa_U}(t)}
		=\gamma^d\Big(\prod_{i=1}^da_i\Big)\Big(\prod_{i=1}^dw_i\Big)^{-\gamma-1}\lambda_U(a_1w_1^{-\gamma}, \dots, a_dw_d^{-\gamma}; \kappa_U)\label{tail density transform}\\
		&=&\lambda_U(a_1w_1^{-\gamma}, \dots, a_dw_d^{-\gamma}; \kappa_U)|J(a_1w_1^{-\gamma}, \dots, a_dw_d^{-\gamma})|,\nonumber
	\end{eqnarray}
	where $V(t)\in \mbox{RV}_{-\gamma}$ and $J(a_1w_1^{-\gamma}, \dots, a_dw_d^{-\gamma})$ is the Jacobian determinant of the homeomorphic transform $y_i=a_iw_i^{-\gamma}$, $1\le i\le d$.
\end{Pro}

\begin{Pro}
	\label{G} 
	\rm  (Gumbel case)
	Assume that marginal densities are right-tail equivalent in the sense of \eqref{equiv}. 
	Assume also that the marginal density $f_i$, $1\le i\le d$, is differentiable and satisfies
	\begin{equation}
		\label{G tail}
		f_i(t+m_i(t)x)\sim f_i(t)e^{-x}\sim a_if_1(t)e^{-x}, ~x\ge 0,~~\mbox{as}~t\to \infty, 
	\end{equation}
	for some self-neglecting function $m_i(\cdot)$ whose derivative converges to $0$; i.e., $m_i(t+m_i(t)x)\sim m_i(t)$, as $t\to \infty$, $1\le i\le d$. If the limit \eqref{U tail} with tail order $\kappa_U$ holds locally uniformly for all ${\boldsymbol w}\in [0,\infty)^d\backslash \{\boldsymbol 0\}$, then the tail density $\lambda(\cdot)$ of $F$ can be expressed in terms of the upper tail density $\lambda_U(\cdot)$ of $C$ as follows: for any ${\boldsymbol w}=(w_1, \dots, w_d)\in \mathbb{R}^d_+\backslash \{\boldsymbol 0\}$, 
	\begin{eqnarray}
		\lambda({\boldsymbol w}) &:=&\lim_{t\to \infty}\frac{f(t+m(t)w_1, \dots, t+m(t)w_d)}{m^{-d}(t)V^{\kappa_U}(t)}\label{tail density transform G}\\
		&=&\Big(\prod_{i=1}^da_i\Big)e^{-\sum_{i=1}^dw_i}\lambda_U(a_1e^{-w_1}, \dots, a_de^{-w_d}; \kappa_U)\nonumber\\
		&=&\lambda_U(a_1e^{-w_1}, \dots, a_de^{-w_d}; \kappa_U)|J(a_1e^{-w_1}, \dots, a_de^{-w_d})|,\nonumber
	\end{eqnarray}
	where $V(t)$ is a function satisfying that $V(t+m(t)x)\sim V(t)e^{-x}$ as $t\to \infty$ for some self-neglecting function $m(t)$, and 
	$J(a_1e^{-w_1}, \dots, a_de^{-w_d})$ is the Jacobian determinant of the homeomorphic transform $y_i=a_ie^{-w_i}$, $1\le i\le d$. 
\end{Pro}

\begin{Rem}
	\begin{enumerate}
		\item These results were first established in \cite{LH14} for the case that $a_i=1$ for $1\le i\le d$. Propositions \ref{F} and \ref{G} under a more general assumption \eqref{equiv} can be proved using the same idea presented in Propositions 2.1 and 2.2 of \cite{LH14} (that is, using line-by-line rewrites of the proofs of Propositions 2.1 and 2.2 of \cite{LH14}). 
		\item The self-neglecting function $m(t)$ in \eqref{tail density transform G} can be taken as $m(t) =m_1(t) = \overline{F}_1(t)/f_1(t)$, $t\ge 0$. The function $V(t)$ in both \eqref{tail density transform} and \eqref{tail density transform G} can be taken as $\Fbar_1(t)
		[\ell(\Fbar_1(t))]^{1/\kappa_U}$, for some function  $\ell(\cdot)$ that is slowly varying at $0$.
		\item  Obviously the tail density $\lambda(\cdot)$ in \eqref{tail density transform} is homogeneous of order $-\gamma \kappa_U-d$. 
	 In contrast, the tail density $\lambda(\cdot)$ in \eqref{tail density transform G} is not homogeneous, and it satisfies the following stability condition:
		\[\lambda(\ww + z{\boldsymbol 1})=\lambda(\ww)e^{-z\kappa_U}, \ {\boldsymbol w}>{\boldsymbol 0}, 
		\]
		for any $z>0$. 
	\end{enumerate}
\end{Rem}

\section{Tail densities of skew-elliptical copulas}
\label{sec-theory}

Propositions \ref{F} and \ref{G}
in the preceding section
can be used to establish the results of skew-elliptical copulas that involve density generators of skew-elliptical distributions.
In this section,
we derive the tail density of \eqref{density Y} at which $f_{\YY}(\yy)$ decays to zero as $\yy$ moves to infinity, 
and  
obtain the scale-invariant tail densities of the skew-elliptical copula densities of $SE_{d}({  \mmu}, \Sigma, g_{d+1}, \ddelta)$ under various assumptions on density generator $g_{d+1}(\cdot)$. 
In contrast to the cumulative-distribution approach used in the multivariate extreme value theory, our tail density approach is especially tractable for the skew-elliptical distribution \eqref{density Y} because the elliptical density $f_{\XX}(\cdot; \mmu, \Sigma)$ has an explicit expression and the cumulative distribution $F(\cdot\ ; g_{q(\yy)})$ is only one-dimensional. 
 It must be mentioned that copula tail densities depend on tail orders and tail behaviors of slowly varying functions in \eqref{tail density transform} or in \eqref{tail density transform G}, extracted from skew-elliptical distributions.

\begin{Pro}
	\label{RV tail}
	\rm
	Let $\YY$ denote a random vector with  skew-elliptical density $f_{\YY}(\cdot)$ specified by \eqref{density Y} or \eqref{density Y explicit} having univariate marginal density  distributions $f_1, \dots, f_d$. If the density generator $g_{d+1}\in \mbox{RV}_{-\alpha}$, $\alpha>(d+1)/2$, then for any $\ww>\Zero$, 
		\begin{equation}
		\label{RV tail eq}
		f_{\YY}(t\ww)\sim 2|\Sigma|^{-1/2} t\,g_{d+1}(t^2)\int_{-\infty}^{\ww\ttheta^\top}(r^2+\ww\Sigma^{-1}\ww^\top)^{-\alpha}{\rm d}r, ~\mbox{as}~t\to \infty.
		\end{equation}
\end{Pro}

The regular variation of the density generator yields the asymptotics \eqref{RV tail eq}, which provides the decay rate of $f_{\YY}(t\ww)$ as $t$ goes to infinity and the corresponding copula tail density, as are showed in the next theorem. The proofs of both Proposition \ref{RV tail} and Theorem \ref{RV tail density} are detailed in the Appendix.

\begin{The}
	\label{RV tail density}
	\rm
	Let $\YY$ denote a random vector having the  skew-elliptical density $f_{\YY}(\cdot)$ specified by \eqref{density Y} or \eqref{density Y explicit} with $\Sigma=(\sigma_{ij})$, where $\sigma_{ii}=1$, $1\le i\le d$, and univariate marginal density  distributions $f_1, \dots, f_d$. Assume that the density generator $g_{d+1}\in \mbox{RV}_{-\alpha}$, $\alpha>(d+1)/2$. 
	\begin{enumerate}
		\item $f_1, \dots, f_d$ are right-tail equivalent in the sense of \eqref{equiv}. 
		\item For any $\ww>\Zero$, 
		\begin{equation}
		\label{Elliptical copula tail density}
	\lambda_U(\ww; 1)=\lambda(a_1^{1/\gamma}w_1^{-1/\gamma}, \dots, a_d^{1/\gamma}w_d^{-1/\gamma})|J(a_1^{1/\gamma}w_1^{-1/\gamma}, \dots, a_d^{1/\gamma}w_d^{-1/\gamma})|,
		\end{equation}
		where $\gamma=2\alpha-d-1>0 $, and 
		\begin{equation}
		\label{Elliptical copula tail density2}
			\lambda({\boldsymbol w})=2K|\Sigma|^{-1/2} \int_{-\infty}^{\ww\ttheta^\top}(r^2+\ww\Sigma^{-1}\ww^\top)^{-\alpha}{\rm d}r,\ \ a_i=\frac{\int_{-\infty}^{\bar\theta_i}(r^2+1)^{-\alpha}{\rm d}r}{\int_{-\infty}^{\bar\theta_1}(r^2+1)^{-\alpha}{\rm d}r},  
		\end{equation}
		where $\bar\theta_i= \delta_i/(1-\delta_i^2)^{1/2}$, for $1\le i\le d$,
{$\ttheta$ is defined in \eqref{eq-ttheta}}, and $K$ is a positive constant. 
	\end{enumerate}
\end{The}

The constant $K$ in \eqref{Elliptical copula tail density2} exhibits flexibility of the limiting decay rate expression of tail density $\lambda(\cdot)$, and can be determined in practice (see Section \ref{t example}). 
Note that the lower tail density of the skew-elliptical copula can be obtained from the upper tail density of $[-\XX | X_0>0]$. The following corollary is immediate. 

\begin{Cor}
	\rm\label{RV lower tail}
	Let $\YY$ denote a random vector having the  skew-elliptical density $f_{\YY}(\cdot)$ specified by \eqref{density Y} or \eqref{density Y explicit} with $\Sigma=(\sigma_{ij})$, where $\sigma_{ii}=1$, $1\le i\le d$, and univariate marginal density  distributions $f_1, \dots, f_d$. Assume that the density generator $g_{d+1}\in \mbox{RV}_{-\alpha}$, $\alpha>(d+1)/2$. 
	\begin{enumerate}
		\item $f_1, \dots, f_d$ are left-tail equivalent in the sense of 
		\begin{equation}
		\label{equiv lower}
		f_i(t)\sim a_if_1(t),~\mbox{as}~t\to -\infty,\quad a_i=\frac{\int_{-\infty}^{-\bar\theta_i}(r^2+1)^{-\alpha}{\rm d}r}{\int_{-\infty}^{-\bar\theta_1}(r^2+1)^{-\alpha}{\rm d}r},  ~~1\le i\le d, 
		\end{equation}
		where $\bar\theta_i=\delta_i/(1-\delta_i^2)^{1/2}$, $1\le i\le d$. 
		\item For any $\ww>\Zero$, {$\lambda_L(\ww; 1)$ has the same form as the right-hand side of
\eqref{Elliptical copula tail density} with $a_i$ from \eqref{equiv lower},}
		$\gamma=2\alpha-d-1>0 $, and 
		\begin{equation}
		\label{Elliptical copula tail density2 lower}
		\lambda({\boldsymbol w})=2K|\Sigma|^{-1/2} \int_{-\infty}^{-\ww\ttheta^\top}(r^2+\ww\Sigma^{-1}\ww^\top)^{-\alpha}{\rm d}r,\ \ 
		\end{equation}
		where 
$\ttheta={\ddelta\Sigma^{-1}}/{(1-\ddelta \Sigma^{-1}\ddelta^\top)^{1/2}}$, 
and $K$ is a positive constant.  
	\end{enumerate}
\end{Cor}

\begin{Rem}
	\label{heavy tail remark}
	\begin{enumerate}
		\item 	If $\ddelta = \Zero$, then  $\ttheta=\Zero$ and thus $f_{\YY}(\cdot)=f_{\XX}(\cdot; \mmu, \Sigma)$ (see \eqref{density Y}). In this case, \eqref{Elliptical copula tail density2} becomes 
		\[\lambda({\boldsymbol w})=2K|\Sigma|^{-1/2} \int_0^{\infty}(r^2+\ww\Sigma^{-1}\ww^\top)^{-\alpha}{\rm d}r,\ \ a_i=\frac{\int_0^{\infty}(r^2+1)^{-\alpha}{\rm d}r}{\int_0^{\infty}(r^2+1)^{-\alpha}{\rm d}r}=1,   \ 1\le i\le d, 
		\]
		which, together with \eqref{Elliptical copula tail density}, present the tail density of the (symmetric) elliptical copula. 
		\item The skewness vector $\ddelta \ge \Zero$ may not imply that $\ttheta\ge \Zero$ since $\ttheta$ is a linearly transformed skewness vector. If, however, $\ttheta\ge \Zero$, then it follows from \eqref{Elliptical copula tail density2} and  \eqref{Elliptical copula tail density2 lower} that the skew-elliptical distribution has a heavier joint upper tail than the joint lower tail. 
	\end{enumerate}
\end{Rem}

Applying \eqref{generator integral} to the generator $g_{d+1}$, we have that 
$ \int_0^\infty r^{(d+1)/2-1}g_{d+1}(r){\rm d}r <\infty, $
which yields the tail estimates: as $t\to \infty$, 
\begin{equation}
\label{generator tail bounds}
g_{d+1}(t)\le \kappa\,t^{-(d+1)/2-\epsilon/2},\ \mbox{or}\ g_{d+1}(t^2)\le \kappa\,t^{-(d+1)-\epsilon}, \ \mbox{for a small $\epsilon>0$},
\end{equation}
where $\kappa$ is a constant. 
If  $g_{d+1}\in \mbox{RV}_{-\alpha}$, $\alpha>(d+1)/2$, as we assume in Theorem \ref{RV tail density}, then the condition \eqref{generator tail bounds} is satisfied. Since regular variation with tail index $\alpha>(d+1)/2$ is stronger than the tail estimate \eqref{generator tail bounds}, this tail estimate is not used in Theorem \ref{RV tail density}, but it is needed in the Gumbel case.

{The following definition is needed for the next proposition.}

\begin{Def}\rm
\label{GMD Gumbel}
A Borel-measurable function $g:[0,\infty)\to \mathbb{R}_+$ is said to be in the $d$-dimensional {\it quadratic max-domain of attraction for the Gumbel distribution} with auxiliary function $m(\cdot)$  if for any $d\times d$ positive-definite matrix $Q$ and for any $\xx\in \mathbb{R}^d_+$, 
	\begin{equation}
\label{Gumbel generator}
g\big([t{\boldsymbol 1}+m(t)\xx]Q[t{\boldsymbol 1}+m(t)\xx]^\top\big)\sim 	g\big(t^2{\boldsymbol 1}Q{\boldsymbol 1}^\top\big)\exp\{-\xx Q {\boldsymbol 1}^\top\},\ \mbox{as}\ t\to \infty, 
\end{equation}
for some self-neglecting function $m(\cdot)$; that is, $m(t+m(t)x)\sim m(t)$, as $t\to \infty$. 
\end{Def}

Note that \eqref{Gumbel generator} essentially means that the function $g(\xx Q\xx^\top)$, $\xx\in \mathbb{R}^d_+$, is in a multivariate max-domain  of attraction for the Gumbel distribution. For example, the density generator of the multivariate normal distribution is given by $g(x)=(2\pi)^{-d/2}\exp\{-x/2\}$, which satisfies \eqref{Gumbel generator}, with auxiliary function $m(t)=t^{-1}$, $t>0$.

\begin{Pro}
	\label{Gumbel tail}
	\rm
	Let $\YY$ denote a random vector with  skew-elliptical density $f_{\YY}(\cdot)$ specified by \eqref{density Y} or \eqref{density Y explicit} having univariate marginal density  distributions $f_1, \dots, f_d$. If the density generator $g_{d+1}$ is in the $(d+1)$-dimensional qua{\rm d}ratic max-domain of attraction for the Gumbel distribution with auxiliary function $m(\cdot)$, then for any $\ww\in \mathbb{R}^d_+$, 
	\begin{equation}
	\label{Gumbel tail eq}
		f_{\YY}\big(t{\boldsymbol 1}+m(t)\ww\big) \sim  
		\left\{\begin{array}{ll}
		2|\Sigma|^{-1/2} g_{d}\big(t^2{\boldsymbol 1}\Sigma^{-1}{\boldsymbol 1}^\top\big)\exp\{-\ww\Sigma^{-1}{\boldsymbol 1}^\top\}, & \mbox{if ${\boldsymbol 1}\ttheta^\top>{0}$} \\
		|\Sigma|^{-1/2} g_{d}\big(t^2{\boldsymbol 1}\Sigma^{-1}{\boldsymbol 1}^\top\big)\exp\{-\ww\Sigma^{-1}{\boldsymbol 1}^\top\}, & \mbox{if ${\boldsymbol 1}\ttheta^\top={0}$}\\
			2|\Sigma|^{-1/2} 
			\int_{-\infty}^{t{\boldsymbol 1}\ttheta^\top}g_{d+1}\big(r^2+t^2{\boldsymbol 1}\Sigma^{-1}{\boldsymbol 1}^\top\big){\rm d}r
			\exp\{-\ww\Sigma^{-1}{\boldsymbol 1}^\top\}, & \mbox{if ${\boldsymbol 1}\ttheta^\top<{0}$}
		\end{array}
		\right.
			\end{equation}
	as $t\to \infty$. 
\end{Pro}

The proof can be found in the Appendix. 
The rapid variation of the density generator yields the asymptotics \eqref{Gumbel tail eq}, which provides the decay rate of  $f_{\YY}(t{\boldsymbol 1}+m(t)\ww)$ as $t$ goes to infinity and the corresponding copula tail density. Note that in contrast to regularly varying decay of $f_{\YY}(\cdot)$ in Theorem \ref{RV tail density}, the decay of $f_{\YY}(\cdot)$ in Proposition \ref{Gumbel tail} is faster because $t{\boldsymbol 1}+m(t)\ww$ moves more slowly, as compared to scaling $t\ww$, to infinity. 

\begin{The}
	\label{Gumbel tail density}
	\rm
	Let $\YY$ denote a random vector having the  skew-elliptical density $f_{\YY}(\cdot)$ specified by \eqref{density Y} or \eqref{density Y explicit} with $\Sigma=(\sigma_{ij})$, where $\sigma_{ii}=1$, $1\le i\le d$, and univariate marginal density  distributions $f_1, \dots, f_d$.
	Assume that the density generator $g_{d+1}$ is in the $(d+1)$-dimensional qua{\rm d}ratic max-domain of attraction for the Gumbel distribution with auxiliary function $m(\cdot)$. Let $\bar\theta_i=\delta_i/(1-\delta_i^2)^{1/2}$, $1\le i\le d$. 
	\begin{enumerate}
		\item If $\bar\theta_i\ge 0$, $1\le i\le d$, then $f_1, \dots, f_d$ are right-tail equivalent in the sense of \eqref{equiv}. 
		\item If  $\bar\theta_i< 0$, $1\le i\le d$, then $f_1, \dots, f_d$ are right-tail equivalent in the sense of \eqref{equiv}, provided that $\bar\theta_1= \dots = \bar\theta_d$. 
		\item In the cases of (1) and (2) above, for any $\ww>\Zero$, 
		\begin{equation}
		\label{Elliptical copula tail density Gumbel}
		\lambda_U(\ww; \kappa_U)=\lambda(-\ln(w_1/a_1), \dots, -\ln(w_d/a_d))|J(-\ln(w_1/a_1), \dots, -\ln(w_d/a_d))|
		\end{equation}
		where $\kappa_U={\boldsymbol 1}\Sigma^{-1}{\boldsymbol 1}^\top>0$ and 
		\begin{equation}
		\lambda({\boldsymbol w})=	\left\{\begin{array}{ll}
		2|\Sigma|^{-1/2}\exp\{-\ww\Sigma^{-1}{\boldsymbol 1}^\top\}, & \mbox{if ${\boldsymbol 1}\ttheta^\top>{0}$ or ${\boldsymbol 1}\ttheta^\top<{0}$,}\\
		|\Sigma|^{-1/2}\exp\{-\ww\Sigma^{-1}{\boldsymbol 1}^\top\}, & \mbox{if ${\boldsymbol 1}\ttheta^\top={0}$,}
		\end{array}
		\right. \label{eq-lambda Gumbel}
		\end{equation}
{with $\ttheta$ defined in \eqref{eq-ttheta}.}
	\end{enumerate}
	\end{The}

The proof is detailed in the Appendix. 
The lower tail density of the skew-elliptical copula can be obtained from the upper tail density of $[-\XX | X_0>0]$, as is shown in the following corollary. 

\begin{Cor}
	\rm\label{Gumbel lower tail case}
	Let $\YY$ denote a random vector having the  skew-elliptical density $f_{\YY}(\cdot)$ specified by \eqref{density Y} or \eqref{density Y explicit} with $\Sigma=(\sigma_{ij})$, where $\sigma_{ii}=1$, $1\le i\le d$, and univariate marginal density  distributions $f_1, \dots, f_d$. 	Assume that the density generator $g_{d+1}$ is in the $(d+1)$-dimensional qua{\rm d}ratic max-domain of attraction for the Gumbel distribution with auxiliary function $m(\cdot)$. 
	\begin{enumerate}
		\item $f_1, \dots, f_d$ are left-tail equivalent in the sense of 
		\begin{equation}
		\label{Gumbel equiv lower}
		f_i(t)\sim a_if_1(t),~\mbox{as}~t\to -\infty,  ~~1\le i\le d, 
		\end{equation}
		if $\bar\theta_i=\delta_i/(1-\delta_i^2)^{1/2}\le 0$, $1\le i\le d$. In addition, $f_1, \dots, f_d$ are left-tail equivalent in the sense of \eqref{Gumbel equiv lower} if $\bar\theta_1= \dots = \bar\theta_d>0$. 
		\item For any $\ww>\Zero$, 
{$\lambda_L(\ww; \kappa_L)$ has the same form as the right-hand side of
\eqref{Elliptical copula tail density Gumbel} with}
 {$a_i$ from \eqref{Gumbel equiv lower},}
		$\kappa_L={\boldsymbol 1}\Sigma^{-1}{\boldsymbol 1}^\top>0$, and 
{$\lambda({\boldsymbol w})$ is the same as \eqref{eq-lambda Gumbel}.}
	\end{enumerate}
\end{Cor}

\begin{Rem}
	\label{Gumbel rem}
	\begin{enumerate}
		\item The conditions on $\bar{\theta}_i$, $1\le i\le d$, such as having the same sign, are needed in Theorem \ref{Gumbel tail density} (Corollary \ref{Gumbel lower tail case}) to ensure the univariate margins are right-tail (left-tail) equivalent in deriving tail densities for multivariate maximums (minimums) in the upper orthant (lower orthant) of $\mathbb{R}^d$. In deriving tail densities for multivariate extremes in another orthant of $\mathbb{R}^d$, the marginal skewness parameters $\bar{\theta}_i$ have different signs to ensure tail equivalence of univariate margins in that orthant. This problem is equivalent to analyzing multivariate extremes in the upper or lower orthant of $\mathbb{R}^d$ after orthogonal transforms and the corresponding tail density can then be derived from  Theorem \ref{Gumbel tail density} or Corollary \ref{Gumbel lower tail case}. 
		\item  If $\ddelta = \Zero$, then  $\ttheta=\Zero$ and thus $f_{\YY}(\cdot)=f_{\XX}(\cdot; \mmu, \Sigma)$ (see \eqref{density Y}). In this case, \eqref{Elliptical copula tail density Gumbel} becomes that
		\[
		\lambda_U(\ww; \kappa_U)=\lambda(-\ln(w_1), \dots, -\ln(w_d))|J(-\ln(w_1), \dots, -\ln(w_d))|
		\]
		where $\kappa_U={\boldsymbol 1}\Sigma^{-1}{\boldsymbol 1}^\top>0$ and 
		$\lambda({\boldsymbol w})=
		|\Sigma|^{-1/2}\exp\{-\ww\Sigma^{-1}{\boldsymbol 1}^\top\}$, for any $\ww>\Zero$.
		It is evident from the proof of Proposition \ref{Gumbel tail} (see \eqref{A(t)}) that even the case that $\ttheta^\top\ge  {\boldsymbol 0}$ with only some $\theta_i=0$ is still covered in Theorem \ref{Gumbel tail density} (1) and (3).  
		\item In contrast to the heavy tail case, if $\ttheta\ge \Zero$, then the upper tail density and lower tail density are the same  (see Remark \ref{heavy tail remark} (2)). Furthermore, since
		\[|J(-\ln(w_1/a_1), \dots, -\ln(w_d/a_d))|=|J(-\ln(w_1), \dots, -\ln(w_d))|, \ \forall \ \ww>\Zero,
		\]
		the tail density of a skew-elliptical copula is proportional to that of the symmetric elliptical copula with the same dispersion matrix $\Sigma$, and the proportional constant is given by
		\[2\exp\{-(\ln a_1, \dots, \ln a_d)\Sigma^{-1}{\boldsymbol 1}^\top\}, 
		\]
		that depends only on $\bar{\theta}_i$, $1\le i\le d$. 
	\end{enumerate}
\end{Rem}

		Note that bivariate tail dependence parameters have been obtained in \cite{HL02, Schmidt02} for elliptical distributions and in \cite{FS10, Padoan2011} for skew-t distributions. In contrast, our method focuses on the densities, that have explicit forms for the symmetrical and asymmetrical cases of elliptical distributions.  	 The explicit density expressions allow us to obtain tail densities of skew-elliptical distributions, from which tail dependence parameters are obtained as integrals of the tail densities;
		see Example \ref{t exa}.
		In the heavy tail (regular variation) case, the insights on multivariate margins and generators obtained from \cite{HJ10, Schmidt02, FS10} are all shared by analysis using our method. Importantly, our results shine more light  into the importance of regular variation in the asymmetric setting for skew-elliptical distributions \cite{FS10}; that is, 
			regular variation is crucial in the emergence and presence of noticeably distinct upper and lower tail dependence for skew-elliptical distributions due to the skewness on slowly decayed tails. The techniques of this paper can be used for other families of multivariate densities which have special
			structural representations.

\section{Tail densities of skew-normal copula and skew-t copula}
\label{sec-examples}

{In this section, it is shown how the theory of Section \ref{sec-theory} 
is applied to some skew-normal and skew-t distributions.}
The results obtained in Section \ref{sec-theory} cover the case of symmetric, elliptical copulas, and also cover the situations that the tail densities of skew-elliptical distributions are asymptotically proportional to the tail densities of symmetric, elliptical distributions within a subcone in $\mathbb{R}^d$. Intuitively, if $F((t\yy-\mmu)\ttheta^\top; g_{q(t\yy)})$ for $\yy\in \mathbb{R}^d_+$ in \eqref{density Y} converges, as $t\to \infty$, to a constant that does not depend $\yy$, then deriving the expression of the skew tail density boils down to finding the tail density of a symmetric, elliptical distribution. This especially applies to various versions of skew-normal distributions.

\subsection{Tail densities of skew-normal distributions}

Let $\XX=(Z_0, Z_1, \dots, Z_d)$ be a random vector with $(d+1)$-dimensional  normal distribution having the mean vector of zeros and the positive-definite correlation matrix
\[\mmu^*=(0, 0, \dots, 0),\ 
\Sigma^*
=
\begin{pmatrix}
1 & \ddelta \\
\ddelta^\top & \Sigma
\end{pmatrix}
\]
where $\ddelta=(\delta_1, \dots, \delta_d)$, and the $d\times d$ positive-definite submatrix $\Sigma=(\sigma_{ij})$ with $\sigma_{ii}=1$, $1\le i\le d$. The density generator is given by $g_{d+1}(x)=(2\pi)^{-(d+1)/2}\exp\{-x/2\}$ that satisfies that 
\[g_{d}(x)=2\int_0^\infty g_{d+1}(r^2+x){\rm d}r,\ x\ge 0.
\] 
A multivariate skew-normal distribution is the distribution of $(Y_1, \dots, Y_d)=[(Z_1, \dots, Z_d)|Z_0>0]$; see \cite{AD96} for details on the properties and application of skew-normal distributions.

It follows from \eqref{eq-ttheta} that $g_{q(\yy)}(s) =(2\pi)^{-1/2}\exp\{-s/2\}$, $\bar\theta_i = \delta_i/(1-\delta_i^2)^{1/2}$, $1\le i\le d$, and 
\[\ttheta=(\theta_1, \dots, \theta_d)=\frac{\ddelta\Sigma^{-1}}{(1-\ddelta \Sigma^{-1}\ddelta^\top)^{1/2}}. 
\]
For example, in the bivariate case, for any $-1<\rho<1$, 
\begin{equation}
\label{bi-normal}
\Sigma = 
\begin{pmatrix}
 1 & \rho \\
 \rho & 1
\end{pmatrix}, \ \ 
\Sigma^{-1}=\frac{1}{1-\rho^2}\begin{pmatrix}
1 & -\rho \\
-\rho & 1
\end{pmatrix}.
\end{equation}
The $i$-th univariate marginal density of $Y_i$ is given by
\[f_{Y_i}(y) = 2\phi(y)\Phi(\bar\theta_iy),\ \bar\theta_i=\delta_i/(1-\delta_i^2)^{1/2},\ 1\le i\le d,
\]
where $\phi(\cdot)$ denotes the standard normal density and $\Phi(\cdot)$ denotes the standard normal cumulative distribution function. 
The joint density of $(Y_1, \dots, Y_d)$ is given by
\begin{equation}
\label{skew-normal exa}
f_d(y_1, \dots, y_d)=2\phi_d(y_1, \dots, y_d; \Sigma)\Phi\Big(\sum_{i=1}^d\theta_iy_i\Big)
\end{equation}
where $\phi_d(\cdot; \Sigma)$ is the $d$-dimensional normal density with the mean vector of zeros and correlation matrix $\Sigma$.

Assume that $\ttheta\ge {\boldsymbol 0}$ and $\bar\theta_i\ge 0$, $1\le i\le d$, and we consider the joint upper tail to illustrate our results. 
The marginal densities are tail equivalent in the sense of \eqref{equiv}, with for $i>1$, 
\begin{equation}
\label{a_2}
a_1=1,\ a_i=\lim_{y\to \infty}[f_{Y_i}(y)/f_{Y_1}(y)] = 
\left\{\begin{array}{ll}
1 & \mbox{if $\bar\theta_1=\bar\theta_i=0$ or $\min\{\bar\theta_1, \bar\theta_i\}>0$}\\
1/2 & \mbox{if $\bar\theta_1>0$ and $\bar\theta_i=0$}\\
2 & \mbox{if $\bar\theta_1=0$ and $\bar\theta_i>0$.}
\end{array}
\right.
\end{equation}
Observe also that the density generator $g_{d+1}(x)=(2\pi)^{-(d+1)/2}\exp\{-x/2\}$ satisfies \eqref{Gumbel generator}, with auxiliary function $m(t)=t^{-1}$, $t>0$.  By Theorem \ref{Gumbel tail density}, the upper tail density of $(Y_1, \dots, Y_d)$ is given by 
\by
\lambda(w_1, \dots, w_d)&=&\lim_{t\to \infty}\frac{f_d(t+m(t)w_1,\dots, t+m(t)w_d)}{g_{d}((t,\dots, t)\Sigma^{-1}(t,\dots, t)^\top)}\\
&=&	\left\{\begin{array}{ll}
2|\Sigma|^{-1/2}\exp\{-(w_1,\dots, w_d)\Sigma^{-1}(1,\dots, 1)^\top\}, & \mbox{if $\max\{\theta_1, \dots, \theta_d\}>0$}\\
|\Sigma|^{-1/2}\exp\{-(w_1,\dots, w_d)\Sigma^{-1}(1,\dots, 1)^\top\}, & \mbox{if $\theta_1=\dots =\theta_d=0$.}
\end{array}
\right.
\ey
Note that in the case that $\max\{\theta_1, \dots, \theta_d\}>0$, $\Phi\Big(t\sum_{i=1}^d\theta_iy_i\Big)$ in \eqref{skew-normal exa} converges to 1 as $t\to \infty$ for any $(y_1, \dots, y_d)>\Zero$, and thus $f_d(y_1, \dots, y_d)\sim 2\phi_d(y_1, \dots, y_d; \Sigma)$ when $(y_1, \dots, y_d)$ moves to positive infinity.

In the bivariate case, when $(\delta_1, \delta_2)\ge (0,0)$ and $\rho\ge 0$, then the upper tail density can be written explicitly as 
\begin{equation}
\label{bi normal tail density}
\lambda(w_1, w_2)=
 \left\{\begin{array}{ll}
2(1-\rho^2)^{-1/2}\exp\{-(w_1+w_2)/(1+\rho)\}, & \mbox{if $\theta_1>0$ or $\theta_2>0$}\\
(1-\rho^2)^{-1/2}\exp\{-(w_1+w_2)/(1+\rho)\}, & \mbox{if $\theta_1=\theta_2=0$.}
\end{array}
\right.
\end{equation}
The upper tail order $\kappa_U=(1,1)\Sigma^{-1}(1,1)^\top=2/(1+\rho)$, and 
{from Theorem \ref{Gumbel tail density}},
the upper tail density of the bivariate skew-elliptical copula is given by
\[\lambda_U(w_1, w_2; \kappa_U)
= \left\{\begin{array}{ll}
2(1-\rho^2)^{-1/2}a_2^{-1/(1+\rho)}(w_1w_2)^{-\rho/(1+\rho)}, & \mbox{if $\theta_1>0$ or $\theta_2>0$}\\
(1-\rho^2)^{-1/2}(w_1w_2)^{-\rho/(1+\rho)}, & \mbox{if $\theta_1=\theta_2=0$,}
\end{array}
\right.
\] 
where $a_2$, given by \eqref{a_2}, depends on skewness  parameters $\delta_1, \delta_2$. 
For the joint lower tail,
we need to have left-tail equivalence in applying  Corollary \ref{Gumbel lower tail case}, and thus we need to assume furthermore that $\delta_1=\delta_2$, which implies that $a_1=a_2=1$. In this case, $\kappa_L=2/(1+\rho)$ and the lower tail density is given by 
\[\lambda_L(w_1, w_2; \kappa_U)
= \left\{\begin{array}{ll}
2(1-\rho^2)^{-1/2}(w_1w_2)^{-\rho/(1+\rho)}, & \mbox{if $\theta_1=\theta_2>0$}\\
(1-\rho^2)^{-1/2}(w_1w_2)^{-\rho/(1+\rho)}, & \mbox{if $\theta_1=\theta_2=0$.}
\end{array}
\right.
\]

\subsection{Tail densities of skew-t distributions}
\label{t example}

Let $\XX=(Z_0, Z_1, \dots, Z_d)$ be a random vector with $(d+1)$-dimensional (central) t distribution having the mean vector of zeros and the positive-definite dispersion matrix
\[\mmu^*=(0, 0, \dots, 0),\ 
\Sigma^*
=
\begin{pmatrix}
1 & \ddelta \\
\ddelta^\top & \Sigma
\end{pmatrix}
\]
where $\ddelta=(\delta_1, \dots, \delta_d)$, and the $d\times d$ positive-definite submatrix $\Sigma=(\sigma_{ij})$ with $\sigma_{ii}=1$, $1\le i\le d$. The density generator is given by
\[g_{d+1}(x)=k(\nu, d+1)\big(\nu+x\big)^{-(\nu+d+1)/2}
\]
where $k(\nu, d+1)$ is a constant defined as
\[k(\nu, d+1)=\frac{\Gamma\big((\nu+d+1)/2\big)\nu^{\nu/2}}{\Gamma\big(\nu/2\big)\pi^{(d+1)/2}}, \ \nu>0, d\ge 0, 
\]
 that depends only on the degree of freedom $\nu$ and dimensionality $d+1$.  A multivariate skew-t distribution is the distribution of $(Y_1, \dots, Y_d)=[(Z_1, \dots, Z_d)|Z_0>0]$; see \cite{BD2001} for details on the properties and application of skew-t distributions.

For skewness towards the joint upper tail, we assume again that
 $\bar\theta_i = \delta_i/(1-\delta_i^2)^{1/2}\ge 0$, $1\le i\le d$, and 
\[(\theta_1, \dots, \theta_d)=\frac{\ddelta\Sigma^{-1}}{(1-\ddelta \Sigma^{-1}\ddelta^\top)^{1/2}}\ge \Zero. 
\]
It follows from \eqref{density Y explicit} that for $\yy=(y_1, \dots, y_d)>\Zero$, 
\[f(y_1, \dots, y_d)=2|\Sigma|^{-1/2}k(\nu, d+1)\int_{-\infty}^{\sum_{i=1}^d\theta_iy_i}\Big(\nu+r^2+\yy\Sigma^{-1}\yy^\top\Big)^{-(\nu+d+1)/2}{\rm d}r.
\]
Observe that $g_{d+1}\in \mbox{RV}_{-(\nu+d+1)/2}$, where $\alpha=(\nu+d+1)/2>(d+1)/2$, and thus by Theorem \ref{RV tail density}, the marginal densities of $Y_i=(Z_i|Z_0>0), 1\le i\le d$, are tail equivalent with $a_i$, $1\le i\le d$, specified in \eqref{Elliptical copula tail density2}.

In the case that all the skewness parameters are equal and positive; that is, $\delta_1= \dots = \delta_d=:\delta>0$, we have that $a_i=1$ for $1\le i\le d$. In this case, let $\bar{\theta}= {\delta_1/(1-\delta_1^2)^{1/2}}>0$, and 
the upper tail density of the skew-t distribution is given by
\[\lambda({\boldsymbol w})=2K|\Sigma|^{-1/2} \int_{-\infty}^{\ww\ttheta^\top}(r^2+\ww\Sigma^{-1}\ww^\top)^{-(\nu+d+1)/2}{\rm d}r,
\]
and the constant $K {=K_1}$ is given by 
\begin{eqnarray}
{K_1} &=& \left(2\frac{\Gamma\big((\nu+1)/2\big)\nu^{(\nu-2)/2}}{\Gamma\big(\nu/2\big)\pi^{1/2}}T_{\nu+1}\big(\bar{\theta}(\nu+1)^{1/2}\big)\right)^{-1}
  \nonumber \\
 &=& \frac{\Gamma\big(\nu/2\big)\pi^{1/2}}{2\Gamma\big((\nu+1)/2\big)\nu^{(\nu-2)/2}T_{\nu+1}\big(\bar{\theta}(\nu+1)^{1/2}\big)}
  \label{eq:K1}
\end{eqnarray}
where $T_{\nu+1}(\cdot)$ is the cumulative distribution function of the t distribution with $\nu+1$ degrees of freedom. 
In the symmetrical case, $\ttheta=\Zero$, and the tail density becomes
\[\lambda({\boldsymbol w})=K|\Sigma|^{-1/2} \int_{-\infty}^{+\infty}(r^2+\ww\Sigma^{-1}\ww^\top)^{-(\nu+d+1)/2}{\rm d}r.
\]
The upper tail density of the skew-t copula is given by
\[\lambda_U(\ww; 1)=2K|\Sigma|^{-1/2}\Big(\frac{1}{\nu}\Big)^d\Big(\prod_{i=1}^dw_i\Big)^{-1/\nu-1}\int_{-\infty}^{\ww^{-1/\nu}\ttheta^\top}(r^2+\ww^{-1/\nu}\Sigma^{-1}(\ww^{-1/\nu})^\top)^{-(\nu+d+1)/2}{\rm d}r
\]
where $\ww^{-1/\nu} = (w_1^{-1/\nu}, \dots, w_d^{-1/\nu})>{\boldsymbol 0}$. The lower tail density of the skew-t copula is different only in the upper limit of the integral, and given by
\[\lambda_L(\ww; 1)=2K|\Sigma|^{-1/2}\Big(\frac{1}{\nu}\Big)^d\Big(\prod_{i=1}^dw_i\Big)^{-1/\nu-1}\int_{-\infty}^{-\ww^{-1/\nu}\ttheta^\top}(r^2+\ww^{-1/\nu}\Sigma^{-1}(\ww^{-1/\nu})^\top)^{-(\nu+d+1)/2}{\rm d}r,
\]
with all the same variables as that in the upper tail case.

\begin{Exa}
	\label{t exa}
	\rm
	In contrast to skew-normal distributions, skewness parameters affect tail dependence parameters significantly for the multivariate t distributions due to slower heavy tail decays. 
To illustrate this and also the fact that our results match the results reported in \cite{FS10, Padoan2011}, we consider the calculations of bivariate tail dependence parameters.
For a bivariate skew-t density with $\delta_1,\delta_2$ that are compatible
with $\rho$, $\theta_1=(\delta_1-\rho\delta_2)/D$ and
$\theta_2=(\delta_2-\rho\delta_1)/D$ where
$D=[(1-\rho^2)(1-\rho^2-\delta_1^2-\delta_2^2+2\delta_1\delta_2\rho)]^{1/2}$.
The density of $(Y_1, Y_2)$ can be written as 
 $$f(y_1,y_2)=2t_{2,\nu}(y_1,y_2;\rho)\, T_{\nu+2}\left(
(\theta_1y_1+\theta_2y_2) \Bigl[{\nu+2\over \nu+ (y_1^2+y_2^2-2\rho
	y_1y_2)/(1-\rho^2)}\Bigr]^{1/2} \right),$$ 
where $\Sigma$ is given by \eqref{bi-normal}, $T_{\nu+2}(\cdot)$ is the cumulative distribution of the t distribution with $\nu+2$ degrees of freedom, and the bivariate t density is given by
\begin{eqnarray*}
 t_{2,\nu}(y_1,y_2;\rho)  &=& (1-\rho^2)^{-1/2} {1\over 2\pi}
 \left(1+\frac{y_1^2+y_2^2-2\rho y_1y_2}{(1-\rho^2)\nu}\right)^{-(\nu+2)/2}.  
\end{eqnarray*}
When $i=1, 2$, the $i$-th marginal skew-t density has the estimate at infinity,
\[f_i(y)\sim 2\frac{\Gamma\big((\nu+1)/2\big)\nu^{(\nu-2)/2}}{\Gamma\big(\nu/2\big)\pi^{1/2}}\nu y^{-(\nu+1)}T_{\nu+1}\big(\delta_i(\nu+1)^{1/2}(1-\delta_i^2)^{1/2}\big)
\] 
and the $i$-th marginal survival function has the estimate,
\[\overline{F}_i(y)\sim 2\frac{\Gamma\big((\nu+1)/2\big)\nu^{(\nu-2)/2}}{\Gamma\big(\nu/2\big)\pi^{1/2}} y^{-\nu}T_{\nu+1}\big(\delta_i(\nu+1)^{1/2}(1-\delta_i^2)^{1/2}\big).
\]
The common tail index is $\gamma=\nu$ with $a_1=1$,
 $$a_2:= \lim_{y\to\oo} \frac{f_2(y)}{f_1(y)}=
  { T_{\nu+1}(\delta_2[\nu+1]^{1/2}[1-\delta_2^2]^{1/2}) \over 
    T_{\nu+1}(\delta_1[\nu+1]^{1/2}[1-\delta_1^2]^{1/2}) }.$$

{With $K_1$ given by \eqref{eq:K1},}
the bivariate skew-t density has the following tail density at infinity:
\begin{eqnarray*}
{\lambda(w_1,w_2)} &=& \lim_{t\to\oo} {f(tw_1,tw_2)\over t^{-(\nu+2)}K^{-1} } = {{K_1} (1-\rho^2)^{(\nu+1)/2} \pi^{-1} \nu^{(\nu+2)/2} 
	\bigl(w_1^2+w_2^2-2\rho w_1w_2\bigr)^{-(\nu+2)/2}}  \\
&& \quad \cdot T_{\nu+2}\left(
	(\theta_1 w_1+\theta_2 w_2) \Bigl[{(\nu+2)(1-\rho^2)\over (w_1^2+w_2^2-2\rho
		w_1w_2)}\Bigr]^{1/2} \right),
\end{eqnarray*}
which is consistent with Theorem \ref{RV tail density}. 
Let 
 $$K^*={ (1-\rho^2)^{(\nu+1)/2} \nu^2 \Gamma(\nu/2) \over  
  2\pi^{1/2} \Gamma([\nu+1]/2) \,
  T_{\nu+1}\bigl(\delta_1[\nu+1]^{1/2}[1-\delta_1^2]^{1/2}\bigr) },$$
be a constant with fixed $\delta_1,\delta_2,\rho,\nu$, and 
then, for $r>0$ and $0\le \omega\le \pi/2$, the tail density can be rewritten as 
\begin{eqnarray*}
 \lambda(r\cos\omega,r \sin\omega) &=&
  K^* r^{-(\nu+2)} (1-2\rho \cos\omega \sin\omega)^{-(\nu+2)/2} \\
 &&\quad \cdot T_{\nu+2}\Bigl({(\theta_1\cos\omega+\theta_2\sin\omega)
  \sqrt{(\nu+2)(1-\rho^2)}\over 
  (1-2\rho \cos\omega \sin\omega)^{1/2}}\Bigr) \\
  &:=& h(\omega)\,r^{-(\nu+2)}.
\end{eqnarray*}

The tail density of the bivariate skew-t copula 
has
tail order $\kappa_U=1$ and slowly varying function $\ell_U(\cdot)=1$. With transform to polar coordinates
$s_1=r\cos\omega$,
$s_2=r\sin\omega$, and $\zeta:=\arctan \{ a_2^{1/\nu}\}$, 
the bivariate upper tail dependence parameter (see \cite{JLN10, Joe2014}) can then be calculated as follows,
\begin{eqnarray*} 
 b_U(1,1) &=& \lim_{u\to 0^+}\frac{\mathbb{P}(\overline{F}_1(X_1)>1-u, 
  \overline{F}_2(X_2)>1-u)}{u} =\int_0^1\int_0^1 \lambda_U(w_1, w_2)
  {\rm d}w_1{\rm d}w_2\\
 &=& \int_{a_1^{1/\nu}}^\infty \int_{a_2^{1/\nu}}^\infty
 \lambda(s_1,s_2)  {\rm d}s_1{\rm d}s_2 \\
 &=& \int_0^{\zeta} \int_{a_2^{1/\nu}/\sin\omega}^\infty \!\!
 \lambda(r\cos\omega,r \sin\omega)\,r{\rm d}r{\rm d}\omega 
 + \int_{\zeta}^{\pi/2}\int_{a_1^{1/\nu}/\sin(\pi/2-\omega)}^\infty\!\!\!\!
 \lambda(r\cos\omega,r \sin\omega)\,r{\rm d}r{\rm d}\omega \\
 &=& \int_0^{\zeta} \int_{a_2^{1/\nu}/\sin\omega}^\infty
  h(\omega)\,r^{-(\nu+2)}\,r{\rm d}r{\rm d}\omega
 + \int_{\zeta}^{\pi/2} \int_{a_1^{1/\nu}/\sin(\pi/2-\omega)}^\infty
  h(\omega)\,r^{-(\nu+2)}\,r{\rm d}r{\rm d}\omega \\
 &=& a_2^{-1}\nu^{-1} \int_0^\zeta h(\omega) (\sin\omega)^\nu {\rm d}\omega
 + \nu^{-1} \int_\zeta^{\pi/2} h(\omega) [\sin(\pi/2-\omega)]^\nu {\rm d}\omega.
\end{eqnarray*}
{For the lower tail dependence parameter, use the above expression with
$(\delta_1,\delta_2)$ replaced by $(-\delta_1,-\delta_2)$ and
$(\theta_1,\theta_2)$ replaced by $(-\theta_1,-\theta_2)$.}
The calculations via this formula match the results reported in \cite{FS10, Padoan2011} (see also \cite{Padoan2016}). 
However our expressions for the tail dependence parameters 
have simpler geometric forms and the limits of the integrands are finite. More importantly, our expressions illustrate how the tail dependence parameters are affected by skewness via the transform $\theta_1w_1+\theta_2w_2$ on the tail densities in the heavy tail case. 
The density plots in \cite{Yoshiba2018} of the bivariate skew-t copula with
standard normal margins provide an idea of how the skewness parameter
affects the relative magnitude of the upper and lower tail densities.
\hfill $\Box$
\end{Exa}

\subsection{Tail densities of bivariate skew-normal distributions of \cite{JSA12}}

The method developed in this paper can be applied to the situation where the conditioning variable $X_0$ in \eqref{density Y explicit} can be an elliptically distributed vector. The method can be also applied to the case that the distribution can be a mixture of skew-elliptical distributions. Instead of developing such a formula, we consider the general version of the multivariate skew-normal distributions introduced in \cite{JSA12}, that can have different skewness parameters for different variables. We only illustrate the case where the univariate marginal distributions are the same.

Let $(Z_1, Z_2)$ have a bivariate normal distribution with zero means, unit variances and correlation $\rho\in (-1, 1)$, and let $W$, independent of $(Z_1, Z_2)$,  have a standard normal distribution. Consider 
\[
Y_i= \left\{\begin{array}{ll}
Z_i, & \mbox{if $\eta Z_i>W$}\\
-Z_i, & \mbox{if $\eta Z_i\le W$,}
\end{array}
\right.
\]
where $\eta\ge 0$ is the skewness parameter. Note that $(Y_1, Y_2)$ is a mixture of skew-normal random vectors, with mixing variable $W$ being random.  
The univariate marginal densities are
\[f_{Y_1}(y)=f_{Y_2}(y)=2\phi(y)\,\Phi(\eta y), i=1, 2,
\]
where $\phi(\cdot)$ denotes the standard normal density and $\Phi(\cdot)$ denotes the standard normal cumulative distribution function. The bivariate density of $(Y_1, Y_2)$ is given by
\[f(y_1, y_2)=2\phi_2(y_1,y_2; \rho)\,\Phi(\eta \min\{y_1, y_2\})+2\phi_2(y_1, -y_2; \rho)(\Phi(\eta y_1)+\Phi(\eta y_2)-1)I_{\{\eta y_1+\eta y_2>0\}},
\]
where $\phi_2(\cdot; \rho)$ is the density of the bivariate normal distribution with zero means and correlation $\rho$, and $I_A$ denotes the indicator function of set $A$. If $\eta=0$, then $f(\cdot)$ becomes the symmetric normal density with correlation $\rho$. If $\eta>0$, then use the  auxiliary function $m(t)=t^{-1}$ and we have
\[f(t+y_1/t, t+y_2/t)\sim 2\phi_2(t+y_1/t, t+y_2/t; \rho)+2\phi_2(t+y_1/t, -t-y_2/t; \rho), \ \mbox{as}\ t\to \infty.
\]
That is, for the upper tail case, the tail density of $(Y_1, Y_2)$ is proportional to two normal tail densities. 
It then follows from \eqref{bi normal tail density} that  the upper tail density is given by
\[\lambda(w_1, w_2) =\left\{\begin{array}{ll} \frac{1}{(1-\rho^2)^{1/2}}\exp\{-(w_1+w_2)/(1+\rho)\}, & \mbox{if $\rho\ge 0$}\\
\frac{1}{(1-\rho^2)^{1/2}}\exp\{-(w_1+w_2)/(1-\rho)\}, & \mbox{if $\rho< 0$,}
\end{array}
\right.
\]
That is, 
$$\lambda(w_1, w_2)= \frac{1}{(1-\rho^2)^{1/2}}\exp\{-(w_1+w_2)/(1+|\rho|)\}$$
and the upper tail order $\kappa_U=2/(1+|\rho|)$. 
{By Theorem \ref{Gumbel tail density}},
the upper tail density of the copula of $(Y_1, Y_2)$ is given by
\[\lambda_U(w_1, w_2; \kappa_U)=\frac{1}{(1-\rho^2)^{1/2}}(w_1w_2)^{-\frac{|\rho|}{1+|\rho|}}, 
\]
where $-1<\rho<1$.

\bigskip
\noindent
{\bf Acknowledgments.} The authors would like to thank two referees for their detailed comments that {led to an improved presentation}. This research has been supported by  NSERC Discovery Grant 8698
and NSF grant DMS 1007556.

\appendix

\section{Proofs}
\label{proofs}

\subsection{Proof of Proposition \ref{RV tail}}

\noindent
{\bf Proof.} Note first that $\alpha>(d+1)/2$ is equivalent to that the integral in \eqref{generator integral} is finite for the generator $g_{d+1}$.  Assume, without loss of generality, that $\mmu=\Zero$. 
Consider, from \eqref{density Y explicit}, that for any fixed $\ww>\Zero$,
\be
f_{\YY}(t\ww)&=& 2|\Sigma|^{-1/2}\int_{-\infty}^{t\ww\ttheta^{\top}}g_{d+1}(r^2+t^2\ww\Sigma^{-1}\ww^{\top}){\rm d}r\nonumber\\
&=& 2|\Sigma|^{-1/2}t\int_{-\infty}^{\ww\ttheta^{\top}}g_{d+1}(t^2r^2+t^2\ww\Sigma^{-1}\ww^{\top}){\rm d}r\nonumber\\
&=& 2|\Sigma|^{-1/2}t\int_{-\ww\ttheta^{\top}}^\infty g_{d+1}(t^2r^2+t^2\ww\Sigma^{-1}\ww^{\top}){\rm d}r\nonumber\\
&=&2|\Sigma|^{-1/2}t\,g_{d+1}(t^2)\left[A(t)+B(t)\right]\label{f_Y}
\ee
where 
\by
A(t)&=&\int_{\{r^2+\ww\Sigma^{-1}\ww^{\top}\le  1\}\cap \{-\ww\ttheta^{\top}\le r< \infty\}}\frac{g_{d+1}(t^2(r^2+\ww\Sigma^{-1}\ww^{\top}))}{g_{d+1}(t^2)}{\rm d}r\\
B(t)&=&\int_{\{r^2+\ww\Sigma^{-1}\ww^{\top}> 1\}\cap \{-\ww\ttheta^{\top}\le r< \infty\}}\frac{g_{d+1}(t^2(r^2+\ww\Sigma^{-1}\ww^{\top}))}{g_{d+1}(t^2)}{\rm d}r.
\ey
Since $\{r: r^2+\ww\Sigma^{-1}\ww^{\top}\le  1\}$ is compact and $\ww\Sigma^{-1}\ww^{\top}>0$, the uniform convergence of the integrand in  $A(t)$ for $r^2$ over $(\ww\Sigma^{-1}\ww^{\top}, \infty)$ ensures that we can pass the limit $t\to \infty$ through the integral sign in $A(t)$; that is, 
\by
\lim_{t\to \infty}A(t)&=& \int_{\{r\in \mathbb{R}:\  r^2+\ww\Sigma^{-1}\ww^{\top}\le  1,\  -\ww\ttheta^{\top}\le r< \infty\}}(r^2+\ww\Sigma^{-1}\ww^\top)^{-\alpha}{\rm d}r,
\ey
due to the fact that  $g_{d+1}\in \mbox{RV}_{-\alpha}$. To calculate $\lim_{t\to \infty}B(t)$, we need to use the Potter bounds (see, e.g., \cite{Resnick07}, page 32). For any sufficiently small $\epsilon>0$, there exists $t_0$ such that for $t\ge t_0$ and all $r$ satisfying that $r^2+\ww\Sigma^{-1}\ww^{\top}> 1$, 
\begin{equation}
	\label{Potter}
	(1-\epsilon)(r^2+\ww\Sigma^{-1}\ww^{\top})^{-\alpha-\epsilon}<\frac{g_{d+1}(t^2(r^2+\ww\Sigma^{-1}\ww^{\top}))}{g_{d+1}(t^2)}< (1+\epsilon)(r^2+\ww\Sigma^{-1}\ww^{\top})^{-\alpha+\epsilon}
\end{equation}
Integrating the functions in \eqref{Potter} and taking the limits yield
\begin{equation}
	\label{lower}
	(1-\epsilon)\int_{\{r\in \mathbb{R}:\  r^2+\ww\Sigma^{-1}\ww^{\top}>  1,\  -\ww\ttheta^{\top}\le r< \infty\}}(r^2+\ww\Sigma^{-1}\ww^{\top})^{-\alpha-\epsilon}{\rm d}r
\end{equation}
\[\le \lim_{t\to \infty}\int_{\{r\in \mathbb{R}:\  r^2+\ww\Sigma^{-1}\ww^{\top}>  1,\  -\ww\ttheta^{\top}\le r< \infty\}}\frac{g_{d+1}(t^2(r^2+\ww\Sigma^{-1}\ww^{\top}))}{g_{d+1}(t^2)}{\rm d}r=\lim_{t\to \infty}B(t)
\]

\begin{equation}
	\label{upper}
	\le (1+\epsilon)\int_{\{r\in \mathbb{R}:\  r^2+\ww\Sigma^{-1}\ww^{\top}>  1,\  -\ww\ttheta^{\top}\le r< \infty\}}(r^2+\ww\Sigma^{-1}\ww^{\top})^{-\alpha+\epsilon}{\rm d}r.
\end{equation}
Since $\alpha>1$ and $r^2+\ww\Sigma^{-1}\ww^{\top}>  1$, the Dominated Convergence Theorem implies that the limits of \eqref{lower} and \eqref{upper} as $\epsilon\to 0$ coincide and thus we have
\[\lim_{t\to \infty}B(t) = \int_{\{r\in \mathbb{R}:\  r^2+\ww\Sigma^{-1}\ww^{\top}>  1,\  -\ww\ttheta^{\top}\le r< \infty\}}(r^2+\ww\Sigma^{-1}\ww^{\top})^{-\alpha}{\rm d}r.
\] 
Plugging $\lim_{t\to \infty}A(t)$ and $\lim_{t\to \infty}B(t)$ into \eqref{f_Y} yields \eqref{RV tail eq}. 
\hfill $\Box$

\subsection{Proof of Theorem \ref{RV tail density}}

\noindent
{\bf Proof.} (1) Let $\Sigma = (\sigma_{ij})$, where $\sigma_{ii}=1$, $1\le i\le d$, and $\ddelta = (\delta_1, \dots, \delta_d)$. Thus,  $\bar\theta_i=\delta_i\sigma_{ii}^{-1}/(1-\delta_i^2\sigma_{ii}^{-1})^{1/2}=\delta_i/(1-\delta_i^2)^{1/2}$, $1\le i\le d$. It follows from \eqref{RV tail eq} that for $1\le i\le d$, 
\begin{equation}
	\label{margin RV}
	f_i(t) \sim 2t\,g_{2}(t^2)\int_{-\infty}^{\bar\theta_i}(r^2+1)^{-\alpha}{\rm d}r, ~\mbox{as}~t\to \infty.
\end{equation}
Hence, as $t\to \infty$,
\[\frac{f_i(t)}{f_1(t)}\to \frac{\int_{-\infty}^{\bar\theta_i}(r^2+1)^{-\alpha}{\rm d}r}{\int_{-\infty}^{\bar\theta_1}(r^2+1)^{-\alpha}{\rm d}r}=:a_i, \ 1\le i\le d. 
\]
That is, $f_1, \dots, f_d$ are right-tail equivalent in the sense of \eqref{equiv}. 

(2) We first claim that $g_{d+1}\in \mbox{RV}_{-\alpha}$ implies that $g_{d}\in \mbox{RV}_{-\alpha+1/2}$, where $g_{d}$ is the density generator of a $d$-dimensional margin of $E_{d+1}({  \mmu}^*, \Sigma^*, g_{d+1})$. This follows from the Tauberian theorem on Stieltjes transforms of regular variation (see Section 1.7 in \cite{BGT1987}) and rewriting \eqref{g} as follows, 
\[g_{d}(s)= \int_s^\infty (r-s)^{-1/2}g_{d+1}(r){\rm d}r=2\int_0^\infty g_{d+1}(r+s){\rm d}r^{1/2}.
\]
Using this property $d-1$ times, we have that $g_{2}\in \mbox{RV}_{-\alpha+(d-1)/2}$. Note that if $d$ is odd, then $g_{2}\in \mbox{RV}_{-\alpha+(d-1)/2}$ also follows directly from Karamata's theorem (see page 25 in \cite{Resnick07}) and the recursive expression obtained in page 37 in \cite{FKN90}.

It follows from \eqref{margin RV} that $f_i\in \mbox{RV}_{-2\alpha+d}$, and this and  Karamata's theorem imply that $\overline{F}_i\in \mbox{RV}_{-2\alpha+d+1}$. Compare \eqref{RV tail eq} to \eqref{tail density transform}, we have
\[1-2\alpha = -d+(-2\alpha+d+1)\kappa_U,
\]
which implies that $\kappa_U=1$. Let $V(t):=t^{d+1}g_{d+1}(t^2)/K\in \mbox{RV}_{-\gamma}$, for some constant $K>0$, where $\gamma=2\alpha-d-1>0$, and then by Proposition \ref{RV tail}, the tail density is given by 
\[\lambda({\boldsymbol w})=\lim_{t\to \infty}\frac{f(t{\boldsymbol w})}{t^{-d}V^{\kappa_U}(t)}= 2K|\Sigma|^{-1/2}\int_{-\infty}^{\ww\ttheta^\top}(r^2+\ww\Sigma^{-1}\ww^\top)^{-\alpha}{\rm d}r, 
\]
and by \eqref{tail density transform}, the tail density of the skew-elliptical copula is given by
\[
\lambda_U(\ww; 1)=\lambda(a_1^{1/\gamma}w_1^{-1/\gamma}, \dots, a_d^{1/\gamma}w_d^{-1/\gamma})|J(a_1^{1/\gamma}w_1^{-1/\gamma}, \dots, a_d^{1/\gamma}w_d^{-1/\gamma})|
\]
for any ${\boldsymbol w}=(w_1, \dots, w_d)>{\boldsymbol 0}$. 
\hfill $\Box$

\subsection{Proof of Proposition \ref{Gumbel tail}}

\noindent
{\bf Proof.} Assume, without loss of generality, that $\mmu=\Zero$. 
Consider, from \eqref{density Y explicit}, that for any fixed $\ww>\Zero$,
\be
f_{\YY}\big(t{\boldsymbol 1}+m(t)\ww\big)&=& 2|\Sigma|^{-1/2}\int_{-\infty}^{[t{\boldsymbol 1}+m(t)\ww]\ttheta^{\top}}g_{d+1}\big(r^2+[t{\boldsymbol 1}+m(t)\ww]\Sigma^{-1}[t{\boldsymbol 1}+m(t)\ww]^{\top}\big){\rm d}r\nonumber\\
&=& 2|\Sigma|^{-1/2}(A(t)+B(t)), 
\label{f_Y Gumbel}
\ee
where 
\[A(t)=\int_{-\infty}^{t{\boldsymbol 1}\ttheta^{\top}}g_{d+1}\big(r^2+[t{\boldsymbol 1}+m(t)\ww]\Sigma^{-1}[t{\boldsymbol 1}+m(t)\ww]^{\top}\big){\rm d}r
\]
\[B(t)=\int_{t{\boldsymbol 1}\ttheta^{\top}}^{[t{\boldsymbol 1}+m(t)\ww]\ttheta^{\top}}g_{d+1}\big(r^2+[t{\boldsymbol 1}+m(t)\ww]\Sigma^{-1}[t{\boldsymbol 1}+m(t)\ww]^{\top}\big){\rm d}r.
\]

(1) Consider the case that ${\boldsymbol 1}\ttheta^\top>0$. Since 
\[A(t)\ge \int_{0}^{t{\boldsymbol 1}\ttheta^{\top}}g_{d+1}\big(r^2+[t{\boldsymbol 1}+m(t)\ww]\Sigma^{-1}[t{\boldsymbol 1}+m(t)\ww]^{\top}\big){\rm d}r, \ \mbox{and}\ \lim_{t\to \infty}\frac{m(t)}{t}= 0, 
\]
the mean-value theorem, together with the tail decreasing property of $g_{d+1}$, imply that
\begin{equation}
	\label{B(t) estimate}
	\frac{|B(t)|}{A(t)}\le \frac{m(t)|\ww\ttheta^\top|}{t{\boldsymbol 1}\ttheta^{\top}}\to 0, \ \mbox{as}\ t\to \infty.
\end{equation}
That is, $f_{\YY}\big(t{\boldsymbol 1}+m(t)\ww\big)\sim 2|\Sigma|^{-1/2}A(t)$, as $t\to \infty$. With substitution $r=ts$, we have
\by
A(t) &=& t\,\int_{-\infty}^{{\boldsymbol 1}\ttheta^{\top}}g_{d+1}\big(t^2s^2+[t{\boldsymbol 1}+m(t)\ww]\Sigma^{-1}[t{\boldsymbol 1}+m(t)\ww]^{\top}\big){\rm d}s\\
&=& t\,\int_{-\infty}^{{\boldsymbol 1}\ttheta^{\top}}g_{d+1}\big([t+m(t)\times 0, t{\boldsymbol 1}+m(t)\ww]Q_s[t+m(t)\times 0, t{\boldsymbol 1}+m(t)\ww]^{\top}\big){\rm d}s
\ey 
where 
\begin{equation}
	\label{g recursion 1}
	Q_s= 
	\begin{pmatrix}
		s^2 & {\boldsymbol 0}\\
		{\boldsymbol 0}^\top & \Sigma^{-1}
	\end{pmatrix}, 
\end{equation}
and ${\boldsymbol 0}$ is the $d$-dimensional row vector of zeros. Since $g_{d+1}$ is in the $(d+1)$-dimensional qua{\rm d}ratic max-domain of attraction for the Gumbel distribution with auxiliary function $m(\cdot)$, we have for any small $\epsilon>0$, when $t$ is sufficiently large, 
\begin{eqnarray}
	&&g_{d+1}\big([t, t{\boldsymbol 1}]Q_s[t, t{\boldsymbol 1}]^\top\big) \exp\{-[0,\ww] Q_s[1, {\boldsymbol 1}]^\top\}(1-\epsilon) \nonumber \\
	&\le& 
	g_{d+1}\big([t+m(t)\times 0, t{\boldsymbol 1}+m(t)\ww]Q_s[t+m(t)\times 0, t{\boldsymbol 1}+m(t)\ww]^{\top}\big) \label{g recursion} \\
	&\le& g_{d+1}\big([t, t{\boldsymbol 1}]Q_s[t, t{\boldsymbol 1}]^\top\big) \exp\{-[0,\ww] Q_s[1, {\boldsymbol 1}]^\top\}(1+\epsilon),\ \mbox{for any}\ s\in \mathbb{R}.  \nonumber
\end{eqnarray}
Observe that $\exp\{-[0,\ww] Q_s[1, {\boldsymbol 1}]^\top\}=\exp\{-\ww\Sigma^{-1}{\boldsymbol 1}^\top\}$. Taking integrals on expressions in above inequalities and letting $\epsilon \to 0$ lead to 
\be
A(t) &\sim &t\,\int_{-\infty}^{{\boldsymbol 1}\ttheta^{\top}}g_{d+1}\big([t, t{\boldsymbol 1}]Q_s[t, t{\boldsymbol 1}]^\top\big)\exp\{-\ww\Sigma^{-1}{\boldsymbol 1}^\top\} {\rm d}s\nonumber\\
&= &
t\,\exp\{-\ww\Sigma^{-1}{\boldsymbol 1}^\top\}\int_{-\infty}^{{\boldsymbol 1}\ttheta^{\top}}g_{d+1}\big(t^2s^2+t{\boldsymbol 1}\Sigma^{-1}t{\boldsymbol 1}^\top\big){\rm d}s\nonumber\\
&=& \exp\{-\ww\Sigma^{-1}{\boldsymbol 1}^\top\}\int_{-\infty}^{t{\boldsymbol 1}\ttheta^{\top}}g_{d+1}\big(r^2+t{\boldsymbol 1}\Sigma^{-1}t{\boldsymbol 1}^\top\big){\rm d}r\label{A(t)}\\
&=& \exp\{-\ww\Sigma^{-1}{\boldsymbol 1}^\top\}\left[\int_{-\infty}^{0}g_{d+1}\big(r^2+t{\boldsymbol 1}\Sigma^{-1}t{\boldsymbol 1}^\top\big){\rm d}r + \int_{0}^{t{\boldsymbol 1}\ttheta^{\top}}g_{d+1}\big(r^2+t{\boldsymbol 1}\Sigma^{-1}t{\boldsymbol 1}^\top\big){\rm d}r\nonumber
\right].
\ee
Since $g_{d+1}$ is integrable  (see \eqref{generator tail bounds}), we have
\[\frac{\int_{t{\boldsymbol 1}\ttheta^{\top}}^{\infty}g_{d+1}\big(r^2+t{\boldsymbol 1}\Sigma^{-1}t{\boldsymbol 1}^\top\big){\rm d}r}{\int_0^{\infty}g_{d+1}\big(r^2+t{\boldsymbol 1}\Sigma^{-1}t{\boldsymbol 1}^\top\big){\rm d}r}\to 0, \ \mbox{as}\ t\to \infty.
\]
This implies that 
\begin{equation}
	\label{Gumbel tail 3}
	\int_0^{t{\boldsymbol 1}\ttheta^{\top}}g_{d+1}\big(r^2+t{\boldsymbol 1}\Sigma^{-1}t{\boldsymbol 1}^\top\big){\rm d}r\sim 
	\int_0^{\infty}g_{d+1}\big(r^2+t{\boldsymbol 1}\Sigma^{-1}t{\boldsymbol 1}^\top\big){\rm d}r.
\end{equation}
Since $g_{d+1}\big(r^2+t{\boldsymbol 1}\Sigma^{-1}t{\boldsymbol 1}^\top\big)$ is symmetric in $r$, we have
\by
A(t)&\sim & 2\exp\{-\ww\Sigma^{-1}{\boldsymbol 1}^\top\}\int_0^{\infty}g_{d+1}\big(r^2+t{\boldsymbol 1}\Sigma^{-1}t{\boldsymbol 1}^\top\big){\rm d}r \\
&=& \exp\{-\ww\Sigma^{-1}{\boldsymbol 1}^\top\}g_{d}\big(t{\boldsymbol 1}\Sigma^{-1}t{\boldsymbol 1}^\top\big),
\ey
where the equality follows from \eqref{g}. Therefore, 
\[f_{\YY}\big(t{\boldsymbol 1}+m(t)\ww\big)\sim 2|\Sigma|^{-1/2}\exp\{-\ww\Sigma^{-1}{\boldsymbol 1}^\top\}g_{d}\big(t{\boldsymbol 1}\Sigma^{-1}t{\boldsymbol 1}^\top\big), 
\]
as $t\to \infty$.

(2) Consider the case that ${\boldsymbol 1}\ttheta^\top={0}$. It follows immediately from the mean-value theorem, \eqref{generator tail bounds} and \eqref{g} that $|B(t)|/A(t)$ converges to zero as $t\to \infty$.  In this situation, \eqref{f_Y Gumbel} becomes 
\[f_{\YY}\big(t{\boldsymbol 1}+m(t)\ww\big)
\sim 2|\Sigma|^{-1/2}A(t), 
\]
where 
\[A(t)=\int_{-\infty}^{0}g_{d+1}\big(r^2+[t{\boldsymbol 1}+m(t)\ww]\Sigma^{-1}[t{\boldsymbol 1}+m(t)\ww]^{\top}\big){\rm d}r.
\]
It follows from \eqref{A(t)} that 
\[A(t)\sim \exp\{-\ww\Sigma^{-1}{\boldsymbol 1}^\top\}\int_{-\infty}^{0}g_{d+1}\big(r^2+t{\boldsymbol 1}\Sigma^{-1}t{\boldsymbol 1}^\top\big){\rm d}r=\frac{1}{2}\exp\{-\ww\Sigma^{-1}{\boldsymbol 1}^\top\}g_{d}\big(t{\boldsymbol 1}\Sigma^{-1}t{\boldsymbol 1}^\top\big). 
\]
That is,
\[f_{\YY}\big(t{\boldsymbol 1}+m(t)\ww\big)\sim |\Sigma|^{-1/2}\exp\{-\ww\Sigma^{-1}{\boldsymbol 1}^\top\}g_{d}\big(t{\boldsymbol 1}\Sigma^{-1}t{\boldsymbol 1}^\top\big), 
\]
as $t\to \infty$. 

(3) Consider the case that ${\boldsymbol 1}\ttheta^\top< {0}$. Note that in this case, $[t{\boldsymbol 1}+m(t)\ww]\ttheta^{\top}$ can be smaller than $t{\boldsymbol 1}\ttheta^\top$. Since $m(t)/t\to 0$ as $t\to \infty$, the estimate on $B(t)$ in \eqref{B(t) estimate} is still correct. Thus $f_{\YY}\big(t{\boldsymbol 1}+m(t)\ww\big)\sim 2|\Sigma|^{-1/2}A(t)$. Using the same idea as in \eqref{g recursion}, we have
\[f_{\YY}\big(t{\boldsymbol 1}+m(t)\ww\big)\sim 2|\Sigma|^{-1/2}\exp\{-\ww\Sigma^{-1}{\boldsymbol 1}^\top\}\int_{-\infty}^{t{\boldsymbol 1}\ttheta^{\top}}g_{d+1}\big(r^2+t{\boldsymbol 1}\Sigma^{-1}t{\boldsymbol 1}^\top\big){\rm d}r, 
\]
as $t\to \infty$. 
\hfill $\Box$

\subsection{Proof of Theorem \ref{Gumbel tail density}}

\noindent
{\bf Proof.} (1) It follows from \eqref{density Y explicit} with $\mmu = {\boldsymbol 0}$ and \eqref{g} that if $\bar\theta_i=\delta_i/(1-\delta_i^2)^{1/2}$, then
\begin{equation}
	\label{Gumbel tail 4}
	f_{Y_i}(t) = 2\int_{-\infty}^{t\bar\theta_i}g_{2}(r^2+t^2){\rm d}r\sim 4\int_{0}^\infty g_{2}(r^2+t^2){\rm d}r =2g_{1}(t^2), 
\end{equation}
where the tail equivalence follows from the fact that $\int_{0}^{t\bar\theta_i} g_{2}(r^2+t^2){\rm d}r\sim \int_{0}^\infty g_{2}(r^2+t^2){\rm d}r$ for $\bar\theta_i>0$ (see \eqref{Gumbel tail 3}). If $\bar\theta_i=0$, then
\begin{equation}
	\label{Gumbel tail 5}
	f_{Y_i}(t) = 2\int_{-\infty}^{0}g_{2}(r^2+t^2){\rm d}r= 2\int_{0}^\infty g_{2}(r^2+t^2){\rm d}r =g_{1}(t^2). 
\end{equation}
Thus $f_1, \dots, f_d$ are right-tail equivalent in the sense of \eqref{equiv}. If $\bar\theta_1=0$, then for $1\le i\le d$, 
\[a_i=\left\{\begin{array}{ll}
1 & \mbox{if $\bar\theta_i=0$}\\
2 & \mbox{if $\bar\theta_i>0$.}
\end{array}
\right.
\]
If $\bar\theta_1>0$, then for $1\le i\le d$, 
\[a_i=\left\{\begin{array}{ll}
1/2 & \mbox{if $\bar\theta_i=0$}\\
1 & \mbox{if $\bar\theta_i>0$.}
\end{array}
\right.
\]

(2) In the case that $\bar\theta_i< 0$, $1\le i\le d$, we have
\begin{equation}
	\label{Gumbel tail 6}
	f_{Y_i}(t) = 2\int_{-\infty}^{t\bar\theta_i}g_{2}(r^2+t^2){\rm d}r, \ i=1, \dots, d.
\end{equation}
Using L'Hospital's Rule, $f_{Y_i}(t)/f_{Y_1}(t)\sim g_{2}(t^2(\bar\theta_i^2+1))/g_{2}(t^2(\bar\theta_1^2+1))\to a_i=1$, as $t\to \infty$. 

(3) According to \eqref{Gumbel tail eq}, the marginal density of $Y_i$ satisfies that for $x>0$, 
\[
f_{Y_i}(t+m(t)x)\sim 
\left\{\begin{array}{ll}
2g_{1}(t^2)e^{-x} & \mbox{if $\bar\theta_i>0$}\\
g_{1}(t^2)e^{-x} & \mbox{if $\bar\theta_i=0$}\\
2\int_{-\infty}^{t\bar\theta_i}g_{2}(r^2+t^2){\rm d}re^{-x} & \mbox{if $\bar\theta_i<0$.}
\end{array}
\right.
\]
This, together with \eqref{Gumbel tail 4}, \eqref{Gumbel tail 5} and \eqref{Gumbel tail 6}, imply that for any $1\le i\le d$, 
\[f_{Y_i}(t+m(t)x)\sim f_{Y_i}(t)e^{-x}, \ x>0.
\]
Thus the marginal distribution $F_i$ is in the max-domain of attraction of the Gumbel distribution and the reciprocal hazard rate $\overline{F}_1(t)/f_1(t)\sim \overline{F}_i(t)/f_i(t)$, $1\le i\le d$, can be taken as $m(t)$, $t\ge 0$. 
We now show that the density generator $g_{k}$, $k\le d$,  is in the $k$-dimensional qua{\rm d}ratic max-domain of attraction for the Gumbel distribution with auxiliary function $m(\cdot)$. Taking integrals with respect to $s$ on the functions in \eqref{g recursion} and letting $\epsilon\to 0$ lead to 
\[\int_0^\infty g_{d+1}\big(r^2+[t{\boldsymbol 1}+m(t)\ww]\Sigma^{-1}[t{\boldsymbol 1}+m(t)\ww]^{\top}\big){\rm d}r \]
\begin{equation}
	\label{g recursion 2}
	\sim \exp\{-\ww\Sigma^{-1}{\boldsymbol 1}^\top\}\int_0^\infty g_{d+1}\big(r^2+t{\boldsymbol 1}\Sigma^{-1}t{\boldsymbol 1}^{\top}\big){\rm d}r,\ t\to \infty.
\end{equation}
It then follows from \eqref{g} that for any $d\times d$ positive-definite matrix $\Sigma^{-1}$, 
\begin{equation}
	\label{g tail}
	g_{d}\big([t{\boldsymbol 1}+m(t)\ww]\Sigma^{-1}[t{\boldsymbol 1}+m(t)\ww]^{\top} \big)\sim \exp\{-\ww\Sigma^{-1}{\boldsymbol 1}^\top\}g_{d}\big(t{\boldsymbol 1}\Sigma^{-1}t{\boldsymbol 1}^{\top}\big),\ t\to \infty.
\end{equation}
That is, the density generator $g_{d}$ is in the $d$-dimensional qua{\rm d}ratic max-domain of attraction for the Gumbel distribution with auxiliary function $m(\cdot)$. The repeated uses of \eqref{g recursion 2} and \eqref{g tail} yield the fact that the density generator $g_{k}$, $k\le d$,  is in the $k$-dimensional qua{\rm d}ratic max-domain of attraction for the Gumbel distribution with auxiliary function $m(\cdot)$. In particular,  the function $G(t):=g_{d}\big(t{\boldsymbol 1}\Sigma^{-1}t{\boldsymbol 1}^{\top}\big)$ satisfies that for any $x>0$, 
\by
G(t+m(t)x) &=& g_{d}\big([t{\boldsymbol 1}+m(t)x{\boldsymbol 1}]\Sigma^{-1}[t{\boldsymbol 1}+m(t)x{\boldsymbol 1}]^{\top} \big)\\
&\sim &\exp\{-x{\boldsymbol 1}\Sigma^{-1}{\boldsymbol 1}^\top\}g_{d}\big(t{\boldsymbol 1}\Sigma^{-1}t{\boldsymbol 1}^{\top}\big)=\exp\{-x{\boldsymbol 1}\Sigma^{-1}{\boldsymbol 1}^\top\}G(t), \ \mbox{as}\ t\to \infty. 
\ey
That is, $G^{1/\kappa_U}(t+m(t)x)\sim e^{-x}G^{1/\kappa_U}(t)$, $t\to \infty$, where $\kappa_U={\boldsymbol 1}\Sigma^{-1}{\boldsymbol 1}^\top>0$. Let
\[V(t):= m^{d/\kappa_U}(t)G^{1/\kappa_U}(t), t\ge 0,
\]
and we have $V(t+m(t)x)\sim e^{-x}V(t)$, $t\to \infty$, and $m^{-d}(t)V^{\kappa_U}(t)=g_{d}\big(t{\boldsymbol 1}\Sigma^{-1}t{\boldsymbol 1}^{\top}\big)$. It now follows from Proposition \ref{G} and \eqref{Gumbel tail eq} that
\[
\lambda({\boldsymbol w}) :=\lim_{t\to \infty}\frac{f(t+m(t)w_1, \dots, t+m(t)w_d)}{m^{-d}(t)V^{\kappa_U}(t)}=
\left\{\begin{array}{ll}
2|\Sigma|^{-1/2}\exp\{-\ww\Sigma^{-1}{\boldsymbol 1}^\top\}, & \mbox{if ${\boldsymbol 1}\ttheta^\top>0$}\\
|\Sigma|^{-1/2}\exp\{-\ww\Sigma^{-1}{\boldsymbol 1}^\top\}, & \mbox{if ${\boldsymbol 1}\ttheta^\top={0}$.}
\end{array}
\right.
\]
The tail density of the copula of $\YY$ is then obtained as in \eqref{Elliptical copula tail density Gumbel}.

In the case that ${\boldsymbol 1}\ttheta^\top< {0}$. 
Let $G(t) := \int_{-\infty}^{t{\boldsymbol 1}\ttheta^{\top}}g_{d+1}\big(r^2+t{\boldsymbol 1}\Sigma^{-1}t{\boldsymbol 1}^\top\big){\rm d}r$. Using \eqref{g recursion}, the univariate function $G(t)$ satisfies that 
\[
G(t+m(t)x)\sim \int_{-\infty}^{t{\boldsymbol 1}\ttheta^{\top}}g_{d+1}\big(r^2+[t{\boldsymbol 1}+m(t)x{\boldsymbol 1}]\Sigma^{-1}[t{\boldsymbol 1}+m(t)x{\boldsymbol 1}]^\top\big){\rm d}r=\exp{\{-x{\boldsymbol 1}\Sigma^{-1}{\boldsymbol 1}^\top\}}G(t).
\] 
That is, $G^{1/\kappa_U}(t+m(t)x)\sim e^{-x}G^{1/\kappa_U}(t)$, $t\to \infty$, where $\kappa_U={\boldsymbol 1}\Sigma^{-1}{\boldsymbol 1}^\top>0$. Let
\[V(t):= m^{d/\kappa_U}(t)G^{1/\kappa_U}(t), t\ge 0,
\]
and we have that $V(t+m(t)x)\sim e^{-x}V(t)$, $t\to \infty$, and $m^{-d}(t)V^{\kappa_U}(t)=G(t)$, and \eqref{Elliptical copula tail density Gumbel} follows in this case. 
\hfill $\Box$


\end{document}